\documentclass[11pt]{article}
\usepackage{makeidx}
\usepackage{amsfonts}
\usepackage{amsmath}
\usepackage{latexsym}
\usepackage{amssymb}
\usepackage{amsthm}
\usepackage{amscd}
\usepackage[dvips]{graphicx}
\usepackage[dvips]{color}

\makeatletter
\newfont{\footsc}{cmcsc10 at 8truept}
\newfont{\footbf}{cmbx10 at 8truept}
\newfont{\footrm}{cmr10 at 10truept}
\theoremstyle{definition}
\newtheorem{theorem}{Theorem}[section]
\newtheorem{prop}[theorem]{Proposition}
\newtheorem{lemma}[theorem]{Lemma}
\newtheorem{corollary}[theorem]{Corollary}

\newtheorem{remark}[theorem]{Remark}

\newenvironment{demo}[1]{%
  \trivlist
  \item[\hskip\labelsep
        {\bf #1.}]
}{%
  \endtrivlist
}
\numberwithin{equation}{section}

\newdimen\Squaresize \Squaresize=20pt
\newdimen\thickness \thickness=1pt         
                                                    
\def\Square#1{\hbox{\vrule width \thickness
   \vbox to \Squaresize{\hrule height \thickness\vss                            
      \hbox to \Squaresize{\hss#1\hss}
   \vss\hrule height\thickness} 
\unskip\vrule width \thickness} 
\kern-\thickness}                                                            
                               
\def\vsquare#1{\vbox{\Square{$#1$}}\kern-\thickness}
\def\blank{\omit\hskip\Squaresize}

\def\fibyoung#1{\let\\=\cr              %added  let\\ =\cr 
\vbox{\smallskip\offinterlineskip
\halign{&\vsquare{##}\cr #1}}\,}

\def\borderlessrect#1#2{\hbox{\hskip \thickness
   \vbox to \Squaresize{\vskip \thickness \vss
      \hbox to #2 {\hss #1\hss}
   \vss\vskip\thickness} 
\unskip\hskip \thickness} 
\kern-\thickness}                                                            
                               
\def\vborderlessrect#1#2{\vbox{\borderlessrect{$#1$}{#2}}\kern-\thickness}

\def\borderless#1{\omit\vborderlessrect{#1}{\Squaresize}}
\def\borderlessrc#1#2{\omit\vborderlessrect{#1}{#2}}

\def\msquare#1{\vbox{\hbox{\vrule width \thickness
   \vbox to \Squaresize{\hrule height \thickness
      \hbox to \Squaresize{\hfil{\sevenrm #1}}
   \vfil\hrule height\thickness}
\unskip\vrule width \thickness}
\kern-\thickness}\kern-\thickness}

\def\twosquare#1#2{\vbox{\hbox{\vrule width \thickness %\vskip 1pt
   \vbox to \Squaresize{\hrule height \thickness
      \hbox to \Squaresize{\hfil{\sevenrm #1}}\vss
      \hbox to \Squaresize{\hss{#2}\hss}
   \vfil\hrule height\thickness}
\unskip\vrule width \thickness}
\kern-\thickness}\kern-\thickness}

\def\twoblank#1#2{\vbox{\hbox{
   \vbox to \Squaresize{\vskip 2pt
      \hbox to \Squaresize{\hfil{\sevenrm #1}\ }\vss
      \hbox to \Squaresize{\hss{#2}\hss}
   \vfil}\unskip\kern-\thickness}
}\unskip\kern-\thickness}

\def\young#1{
\def\>{\blank}
\def\<{\borderless}
\def\*{\borderlessrc}
\def\p{\omit\msquare}
\def\t{\omit\twosquare}
\def\b{\omit\twoblank}
\let\\=\cr %added  let\\ =\cr 
\vbox{\smallskip\offinterlineskip
\halign{&\vsquare{##}\cr #1}}}

\newdimen\smsquaresize \smsquaresize=12pt
\newdimen\smthickness \smthickness=.5pt
\font\smcellfont=cmss8 scaled \magstep0

\def\smsquare#1{\hbox{\vrule width \smthickness
   \unskip\vbox to \smsquaresize{\hrule height \smthickness\vss
      \hbox to \smsquaresize{\hss{\smcellfont #1}\hss}
   \vss\hrule height\smthickness} 
\unskip\vrule width \smthickness} 
\kern-\smthickness}

\def\smvsquare#1{\vbox{\smsquare{$#1$}}\kern-\smthickness}
\def\blank{\omit\hskip\smsquaresize}

\def\smborderlessrect#1#2{\hbox{\hskip \smthickness
   \vbox to \smsquaresize{\vskip \smthickness \vss
      \hbox to #2 {\hss #1\hss}
   \vss\vskip\smthickness} 
\unskip\hskip \smthickness} 
\kern-\smthickness}                                                            
                               
\def\smvborderlessrect#1#2{\vbox{\smborderlessrect{$#1$}{#2}}\kern-\smthickness}

\def\smborderless#1{\omit\smvborderlessrect{#1}{\smsquaresize}}
\def\smborderlessrc#1#2{\omit\smvborderlessrect{#1}{#2}}

\def\smyoung#1{
\def\<{\smborderless}
\def\*{\smborderlessrc}
\let\\=\cr %added  let\\ =\cr
\vbox{\smallskip\offinterlineskip
\halign{&\smvsquare{##}\cr #1}}}
\newdimen\vsmsquaresize \vsmsquaresize=10pt
\newdimen\vsmthickness \vsmthickness=.5pt
\font\vsmcellfont=cmsl8 scaled \magstep0
\font\vsmletterfont=cmr6 scaled \magstep0

\def\vsmsquare#1{\hbox{\vrule width \vsmthickness
   \unskip\vbox to \vsmsquaresize{\hrule height \vsmthickness\vss
      \hbox to \vsmsquaresize{\hss{\vsmcellfont #1}\hss}
   \vss\hrule height\vsmthickness} 
\unskip\vrule width \vsmthickness} 
\kern-\vsmthickness}
\def\vsmvsquare#1{\vbox{\vsmsquare{#1}}\kern-\vsmthickness}
\def\vsmblank{\omit\hskip\vsmsquaresize}
\def\vsmborderless#1{\hbox{\hskip \vsmthickness\unskip
   \vbox to \vsmsquaresize{\vss
      \hbox to \vsmsquaresize{\hss{\vsmletterfont #1}\hss}
   \vss} 
\unskip\hskip \vsmthickness} 
\kern-\vsmthickness}                                                            \def\vsmvborderless#1{\vbox{\vsmborderless{#1}}\kern-\vsmthickness}

\def\vsmyoung#1{
\def\>{\vsmblank}
\def\<{\omit\vsmvborderless}
\let\\=\cr %added  let\\ =\cr
\vbox{\smallskip\offinterlineskip
\halign{&\vsmvsquare{##}\cr #1}}}
\def\rdots{\mathinner{\mkern1mu\raise1pt\hbox{.}\mkern2mu\raise4pt\hbox{.}\mkern2mu\raise7pt\hbox{.}\mkern1mu}}

\def\defterm#1{{\sl #1}\/}

\def\Pf{\operatorname{Pf}}

\def\sgn{\operatorname{sgn}}
\def\MOD{{\,\operatorname{mod}\,}}

\def\MOD{\operatorname{mod}}

\def\cc#1{{\ooalign{\hfil\raise-.02ex\hbox{#1}\hfil\crcr\mathhexbox20D}}}

\def\PHI#1#2#3#4{
{}_2\phi_{1}\left(
{{{#1},{#2}}\atop{{#3}}}
;
\,q,{#4}
\right)}

\def\PHII#1#2#3#4#5#6{
{}_3\phi_{2}\left(
{{{#1},{#2},{#3}}\atop{{#4},{#5}}}
;
\,q,{#6}
\right)}

\def\wt{\operatorname{wt}}

\title{The Andrews-Stanley partition function\\
and\\
Al-Salam-Chihara polynomials
}

\author{Masao ISHIKAWA\\
\small Faculty of Education, Tottori University\\[-0.8ex]
\small Koyama, Tottori, Japan\\[-0.8ex]
\small \texttt{ishikawa@fed.tottori-u.ac.jp}
\and
Jiang ZENG\\
\small Institut Camille Jordan\\[-0.8ex]
\small Universit\'e Claude Bernard Lyon I\\[-0.8ex]
\small 43, boulevard du 11 novembre 1918\\[-0.8ex]
\small 69622 Villeurbanne Cedex, France\\[-0.8ex]
\small \texttt{zeng@igd.univ-lyon1.fr}
}

\date{
}

\begin{document}
\maketitle

%\begin{center}
%Department of Mathematics,
%Faculty of Education,
%Tottori University\\
%\medskip
%and\\
%\medskip
%
%Institut Camille Jordan,
%Universit\'e Claude Bernard Lyon I
%\end{center}

%\noindent
%{Keywords:}\quad Andrews-Stanley partition function, Al-Salam-Chihara polynomials,
%generating functions, Schur's Q-functions, minor summation formula of Pfaffians.

\abstract{
For any partition $\lambda$ let
 $\omega(\lambda)$ denote the four parameter weight
\begin{equation*}
\omega(\lambda)=a^{\sum_{i\geq1}\lceil\lambda_{2i-1}/2\rceil}
b^{\sum_{i\geq1}\lfloor\lambda_{2i-1}/2\rfloor}
c^{\sum_{i\geq1}\lceil\lambda_{2i}/2\rceil}
d^{\sum_{i\geq1}\lfloor\lambda_{2i}/2\rfloor},
\end{equation*}
and let $\ell(\lambda)$ be the length of $\lambda$.
We show that the generating function $\sum\omega(\lambda)z^{\ell(\lambda)}$, where
the sum runs over all ordinary (resp.
strict) partitions with parts each $\leq N$,
can be  expressed by the Al-Salam-Chihara polynomials.
As a corollary we derive  G.E.~Andrews' result by specializing some parameters and
C.~Boulet's results by letting $N\to +\infty$.
In the last section we prove a Pfaffian formula for the weighted sum
$\sum \omega(\lambda)z^{\ell(\lambda)}P_{\lambda}(x)$ where $P_{\lambda}(x)$
is Schur's $P$-function and
the sum runs over all strict partitions.
}

\bigbreak

%\bigbreak
%\noindent
%2000 Mathematics Subject Classifications:
%05A17, 05A30, 05E05, 05E35,
%33D15, 33D45

%%%%%%%%%%%%%%%%%%%%%%%%%%%%%%%%%%%%%%%%%%%%%%%%%%%%%%%%%%%%%%%%%%%%%%%%%%%%%
%
% Section 1: Introduction
%
%%%%%%%%%%%%%%%%%%%%%%%%%%%%%%%%%%%%%%%%%%%%%%%%%%%%%%%%%%%%%%%%%%%%%%%%%%%%%

\section{Introduction}\label{sec:intro}
For any integer partition $\lambda$, denote by $\lambda'$ its conjugate and
$\ell(\lambda)$ the number of its parts.
Let ${\cal O}(\lambda)$ denote the number of odd parts of $\lambda$
and $|\lambda|$ the sum of its parts.
R.~Stanley (\cite{Stan}) has shown that
if $t(n)$ denotes the number of partitions $\lambda$ of
$n$ for which ${\cal O}(\lambda)\equiv {\cal O}(\lambda')$ ($\MOD4$),
then
\[
t(n)=\frac12\left(p(n)+f(n)\right),
\]
where $p(n)$ is the total number of partitions of $n$,
and $f(n)$ is defined by
$$
\sum_{n=0}^{\infty}f(n)q^n=\prod_{i\geq1}\frac{(1+q^{2i-1})}
{(1-q^{4i})(1+q^{4i-2})}.
$$
Motivated by Stanley's problem,
G.E.~Andrews~\cite{And} assigned the
weight  $z^{{\cal O}(\lambda)} y^{{\cal O}(\lambda')} q^{|\lambda|}$ to
 each partition $\lambda$ and
computed the corresponding generating function of
all partitions with parts each less than or equal to $N$
(see Corollary~\ref{cor:andrews}). The following more general
weight first appeared in Stanley's paper~\cite{Stan2}.
Let $a$, $b$, $c$ and $d$ be commuting indeterminates.
For each partition $\lambda$,
define the \emph{Andrews-Stanley partition functions} $\omega(\lambda)$ by
\begin{equation}
\label{eq:weight}
\omega(\lambda)=a^{\sum_{i\geq1}\lceil\lambda_{2i-1}/2\rceil}
b^{\sum_{i\geq1}\lfloor\lambda_{2i-1}/2\rfloor}
c^{\sum_{i\geq1}\lceil\lambda_{2i}/2\rceil}
d^{\sum_{i\geq1}\lfloor\lambda_{2i}/2\rfloor},
\end{equation}
where $\lceil x\rceil$ (resp. $\lfloor x\rfloor$) stands for the smallest (resp. largest) integer greater (resp. less) than or equal to $x$
for a given real number $x$. Actually it is more convenient to define the above weight
through the Ferrers diagram of $\lambda$: one fills
the $i$th row of the Ferrers diagram alternatively by $a$ and $b$ (resp. $c$ and $d$)
if $i$ is \emph{odd} (resp. \emph{even}), the weight $w(\lambda)$ is
then equal to the product of all the entries in the diagram.
For example,
if $\lambda=(5,4,4,1)$
then
$\omega(\lambda)$ is the product of the entries in the following diagram for $\lambda$.
\[
\young{
a&b&a&b&a\\
c&d&c&d\\
a&b&a&b\\
c\\
}
\]

In \cite{Bou}
C.~Boulet has obtained results for the generating functions
of all ordinary partitions
and all strict partitions with respect to
the weight \eqref{eq:weight}
(see Corollary~\ref{boulet:strict} and Corollary~\ref{boulet:ordinary}).
On the other hand, A. Sills \cite{Sills} has given a combinatorial proof of Andrews'
 result, which has been further generalized by
  A.~Yee~\cite{Yee} by restricting the sum over  partitions
with parts each $\leq N$ and length $\leq M$.

In this paper we shall generalize Boulet's results
by summing the weight function $\omega(\lambda)z^{\ell(\lambda)}$
 over all the ordinary (resp. strict) partitions with parts each $\leq N$.
It turns out that the corresponding generating functions
 are related to the basic hypergeometric series, namely
the Al-Salam-Chihara polynomials and the associated Al-Salam-Chihara polynomials
(see Theorem~\ref{hyper:strict} and Theorem~\ref{hyper:ordinary}).

This paper can be regarded as a succession of \cite{I},
in which one of the authors gave a Pfaffian formula for the weighted sum
$\displaystyle\sum \omega(\lambda)s_{\lambda}(x)$ of the Schur functions $s_{\lambda}(x)$,
where the sum runs over all ordinary partitions $\lambda$,
and settled an open problem by Richard~Stanley.
Though it is not possible to specialize the Schur functions to $z^{\ell(\lambda)}$,
we show in this paper that this approach still works, i.e.,
we can evaluate the weighted sum
$\sum \omega(\lambda)z^{\ell(\lambda)}$ by using Pfaffians
and minor summation formulas as tools (\cite{IW1}, \cite{IW5}),
but, as an after thought,  we also provide alternative combinatorial proofs.

In the last section we show the weighted sum
$\displaystyle\sum \omega(\mu)z^{\ell(\mu)}P_{\mu}(x)$ of Schur's $P$-functions $P_{\mu}(x)$
(when $z=2$, this equals the weighted sum $\sum \omega(\mu)Q_{\mu}(x)$ of Schur's $Q$-functions $Q_{\mu}(x)$)
can be expressed by a Pfaffian
where $\mu$ runs over all strict partitions (with parts each $\leq N$).

%%%%%%%%%%%%%%%%%%%%%%%%%%%%%%%%%%%%%%%%%%%%%%%%%%%%%%%%%%%%%%%%%%%%%%%%%%%%%
%
% Section 2: Preliminaries
%
%%%%%%%%%%%%%%%%%%%%%%%%%%%%%%%%%%%%%%%%%%%%%%%%%%%%%%%%%%%%%%%%%%%%%%%%%%%%%

\section{Preliminaries}\label{sec:msf}

A $q$-shifted factorial is defined by
\[
(a;q)_0=1,
\qquad
(a;q)_n=(1-a)(1-aq)\cdots(1-aq^{n-1}),
\qquad\qquad n=1,2,\dots.
\]
We also define $(a;q)_{\infty}=\prod_{k=0}^{\infty}(1-aq^{k})$.
Since products of $q$-shifted factorials occur very often,
to simplify them we shall use the compact notations
\begin{align*}
&(a_1,\dots,a_m;q)_{n}
=(a_1;q)_{n}\cdots(a_m;q)_{n},
\\
&(a_1,\dots,a_m;q)_{\infty}
=(a_1;q)_{\infty}\cdots(a_m;q)_{\infty}.
\end{align*}
We define an ${}_{r+1}\phi_{r}$ \defterm{basic hypergeometric series}
by
%-------------------------------------------%
% basic hypergeometric series
%-------------------------------------------%
\begin{align*}
{}_{r+1}\phi_{r}\left(
{{a_1,a_2,\dots,a_{r+1}}\atop{b_1,\dots,b_{r}}}
;
\,q,z
\right)
=\sum_{n=0}^{\infty}
\frac{(a_1,a_2,\dots,a_{r+1};q)_{n}}{(q,b_1,\dots,b_r;q)_{n}}
z^{n}.
\end{align*}

%-------------------------------------------%
% Al-Salam-Chihara polynomials
%-------------------------------------------%
The \defterm{Al-Salam-Chihara polynomial} $Q_{n}(x)=Q_{n}(x;\alpha,\beta|q)$ is,
by definition (cf. \cite[p.80]{KS}),
\begin{align*}
Q_{n}(x;\alpha,\beta|q)
&=\frac{(\alpha\beta;q)_{n}}{\alpha^n}\,
\PHII{q^{-n}}{\alpha u}{\alpha u^{-1}}{\alpha\beta}{0}{q},
\\
&=(\alpha u;q)_{n}u^{-n}\,
\PHI{q^{-n}}{\beta u^{-1}}{\alpha^{-1}q^{-n+1}u^{-1}}{\alpha^{-1}q u},
\\
&=(\beta u^{-1};q)_{n}u^{n}\,
\PHI{q^{-n}}{\alpha u}{\beta^{-1}q^{-n+1}u}{\beta^{-1}q u^{-1}},
\end{align*}
where $x=\frac{u+u^{-1}}2$.
This is a specialization of the Askey-Wilson polynomials (see \cite{GR}),
and satisfies the three-term recurrence relation
\begin{equation}
2xQ_{n}(x)=Q_{n+1}(x)+(\alpha+\beta)q^{n}Q_{n}(x)+(1-q^{n})(1-\alpha\beta q^{n-1})Q_{n-1}(x),
\label{eq:ASC}
\end{equation}
with $Q_{-1}(x)=0$, $Q_{0}(x)=1$.

We also consider a more general recurrence relation:
%-------------------------------------------%
% Assiciated Al-Salam-Chihara Recurrence eq.
%-------------------------------------------%
\begin{equation}
2x {\widetilde Q}_{n}(x)={\widetilde Q}_{n+1}(x)+(\alpha+\beta)tq^{n}{\widetilde Q}_{n}(x)+(1-tq^{n})(1-t\alpha\beta q^{n-1}){\widetilde Q}_{n-1}(x),
\label{eq:AASC}
\end{equation}
which we call the \defterm{associated Al-Salam-Chihara recurrence relation}.
Put
%-------------------------------------------%
% Two linearly independent solutions
%-------------------------------------------%
\begin{align}
&{\widetilde Q}_{n}^{(1)}(x)
=u^{-n}\,
(t\alpha u;q)_n\,
\PHI{t^{-1}q^{-n}}{\beta u^{-1}}{t^{-1}\alpha^{-1}q^{-n+1}u^{-1}}{\alpha^{-1}qu},
\label{eq:Q1}
\\
&{\widetilde Q}_{n}^{(2)}(x)
=u^{n}\,
\frac{(tq;q)_n(t\alpha\beta ;q)_n}{(t\beta uq;q)_n}\,
\PHI{t q^{n+1}}{\alpha^{-1}qu}{t\beta q^{n+1} u}{{\alpha}{u}},
\label{eq:Q2}
%&{\widetilde Q}_{n}^{(2)}(x)
%=u^{-n}\,
%\frac{(tq;q)_n(\alpha\beta t;q)_n}{(\beta tu^{-1}q;q)_n}\,
%{}_2\phi_{1}\!\left(
%{{t q^{n+1},\alpha^{-1}u^{-1}q}\atop{\beta t u^{-1}q^{n+1}}}
%\Big|
%q;\frac{\alpha}{u}
%\right),
%\label{eq:Q2}
%\\
%&{\widetilde Q}_{n}^{(3)}(x)
%=u^{-n}\,
%(\alpha\beta t;q)_n\,
%{}_2\phi_{1}\!\left(
%{{\alpha u^{-1},\beta u^{-1}}\atop{q u^{-2}}}
%\Big|
%q;t q^{n+1}
%\right),
%\label{eq:Q3}
\end{align}
where $x=\frac{u+u^{-1}}2$.
In \cite{IR},
Ismail and Rahman have presented two linearly independent solutions of the associated
Askey-Wilson recurrence equation (see also \cite{GIM,GM}).
By specializing the parameters,
we conclude that
${\widetilde Q}_{n}^{(1)}(x)$ and ${\widetilde Q}_{n}^{(2)}(x)$ are two linearly independent solutions of
the associated Al-Salam-Chihara equation \thetag{\ref{eq:AASC}}
(see \cite[p.203]{GM}).
Here, we use this fact and omit the proof.
%has solutions
%${\widetilde Q}_{n}^{(i)}(u)$ and ${\widetilde Q}_{n}^{(i)}(1/u)$,
%$i=1,2,3$,
%i.e.
%\thetag{\ref{eq:Q1}}, \thetag{\ref{eq:Q2}} and \thetag{\ref{eq:Q3}} satisfy \thetag{\ref{eq:AASC}}.
%and the both are convergent when $|q|$ is sufficiently small.
The series \thetag{\ref{eq:Q1}} and \thetag{\ref{eq:Q2}} are convergent if we assume $|u|<1$ and $|q|<|\alpha|<1$ (see \cite[p.204]{GM}).

%-------------------------------------------%
% Casorati determinant
%-------------------------------------------%
Let
\begin{equation}
W_{n}={\widetilde Q}_{n}^{(1)}(x){\widetilde Q}_{n-1}^{(2)}(x)-{\widetilde Q}_{n-1}^{(1)}(x){\widetilde Q}_{n}^{(2)}(x)
\label{eq:Casorati}
\end{equation}
denote the Casorati determinant of the equation \thetag{\ref{eq:AASC}}.
Since ${\widetilde Q}_{n}^{(1)}(x)$ and ${\widetilde Q}_{n}^{(2)}(x)$ both satisfy
the recurrence equation \thetag{\ref{eq:AASC}},
it is easy to see that $W_{n}$ satisfies the recurrence equation
\[
W_{n+1}=(1-tq^{n})(1-t\alpha\beta q^{n-1})W_{n}.
\]
Using this equation recursively, we obtain
\[
W_{n+1}=(tq,t\alpha\beta;q)_{n} W_{1},
\]
which implies
%-------------------------------------------%
% W_1
%-------------------------------------------%
\[
W_{1}=\frac{\lim_{n\rightarrow\infty}W_{n+1}}{(tq,t\alpha\beta;q)_{\infty}} .
\]
Using \thetag{\ref{eq:Q1}} and \thetag{\ref{eq:Q2}},
we obtain
\begin{align*}
\lim_{n\rightarrow\infty}W_{n+1}
=
\frac{u^{-1}(t\alpha u,tq,t\alpha\beta,\beta u;q)_{\infty}}
{(t\beta u q,\alpha u;q)_{\infty}}
\end{align*}
(for the detail, see \cite{IR}).
Thus we conclude that
\begin{align}
W_{1}
=\frac{u^{-1}(t\alpha u,\beta u;q)_{\infty}}
{(\alpha u,t\beta u q;q)_{\infty}}.
\label{eq:W1}
\end{align}
In the following sections
we need to find a polynomial solution of the recurrence equation \thetag{\ref{eq:AASC}}
which satisfies a given initial condition,
say ${\widetilde Q}_{0}(x)=\widetilde Q_{0}$ and ${\widetilde Q}_{1}(x)=\widetilde Q_{1}$.
Since ${\widetilde Q}^{(1)}_{n}(x)$ and ${\widetilde Q}^{(2)}_{n}(x)$ are linearly independent
solutions of \thetag{\ref{eq:AASC}},
this ${\widetilde Q}_{n}(x)$ can be written as a linear combination of these functions,
say
\[
{\widetilde Q}_{n}(x)=C_1\,{\widetilde Q}^{(1)}_{n}(x)+C_2\,{\widetilde Q}^{(2)}_{n}(x).
\]
If we substitute the initial condition ${\widetilde Q}_{0}(x)=\widetilde Q_{0}$ and ${\widetilde Q}_{1}(x)=\widetilde Q_{1}$ into this equation
and solve the linear equation,
then we obtain
\begin{align*}
&C_1=\frac1{W_{1}}\left\{\widetilde Q_{1}{\widetilde Q}^{(2)}_{0}(x) - \widetilde Q_{0}{\widetilde Q}^{(2)}_{1}(x)\right\},\\
&C_2=\frac1{W_{1}}\left\{\widetilde Q_{0}{\widetilde Q}^{(1)}_{1}(x) - \widetilde Q_{1}{\widetilde Q}^{(1)}_{0}(x)\right\}.
\end{align*}
By \thetag{\ref{eq:W1}},
we obtain
%-------------------------------------------%
% Solution under an initail condition
%-------------------------------------------%
\begin{align}
{\widetilde Q}_{n}(x)&=
\frac
{u(\alpha u,t\beta uq;q)_{\infty}}
{(t\alpha u,\beta u;q)_{\infty}}
\left[
\left\{\widetilde Q_{1}{\widetilde Q}^{(2)}_{0}(x) - \widetilde Q_{0}{\widetilde Q}^{(2)}_{1}(x)\right\}
{\widetilde Q}^{(1)}_{n}(x)
\right.\nonumber\\
&\qquad\qquad\qquad\qquad\qquad\left.
+\left\{\widetilde Q_{0}{\widetilde Q}^{(1)}_{1}(x) - \widetilde Q_{1}{\widetilde Q}^{(1)}_{0}(x)\right\}
{\widetilde Q}^{(2)}_{n}(x)
\right]
\label{eq:Q_n}
\end{align}
with
%-------------------------------------------%
% Constants
%-------------------------------------------%
\begin{align*}
{\widetilde Q}^{(1)}_{0}(x)
&=
\PHI{t^{-1}}{\beta u^{-1}}{t^{-1}\alpha^{-1}u^{-1}q}{\alpha^{-1} u q},
\\
{\widetilde Q}^{(1)}_{1}(x)
&=u^{-1}(1-\alpha tu)\,
\PHI{t^{-1}q^{-1}}{\beta u^{-1}}{t^{-1}\alpha^{-1}u^{-1}}{\alpha^{-1} u q},
\\
{\widetilde Q}^{(2)}_{0}(x)
&=
\PHI{t q}{\alpha^{-1} u q}{t\beta u q}{\alpha u},
\\
{\widetilde Q}^{(2)}_{1}(x)
&=
\frac{u(1-tq)(1-t\alpha\beta)}{(1-t\beta uq)}\,
\PHI{t q^{2}}{\alpha^{-1} u q}{t\beta u q^{2}}{\alpha u}.
\end{align*}
%-------------------------------------------%
% Limit
%-------------------------------------------%
Since
\begin{align*}
\lim_{n\rightarrow\infty}u^{n}\,{\widetilde Q}^{(1)}_{n}(x)
&=\frac{(t\alpha u,\beta u;q)_{\infty}}{(u^2;q)_{\infty}},
\\
\lim_{n\rightarrow\infty}u^{n}\,{\widetilde Q}^{(2)}_{n}(x)
&=0,
\end{align*}
if we take the limit $\lim\limits_{n\to\infty}u^{n}{\widetilde Q}_{n}(x)$,
then we have
\begin{equation}
\lim\limits_{n\to\infty}u^{n}{\widetilde Q}_{n}(x)
=\frac
{u(t\beta uq,\alpha u;q)_{\infty}}
{(u^2;q)_{\infty}}
\left\{\widetilde Q_{1}{\widetilde Q}^{(2)}_{0}(x) - \widetilde Q_{0}{\widetilde Q}^{(2)}_{1}(x)\right\}.
\label{eq:limit}
\end{equation}

In the later half of this section,
we briefly recall our tools,
i.e. partitions and Pfaffians.
We follow the notation in \cite{Ma} concerning partitions and the symmetric functions.
For more
information about the general theory of determinants and Pfaffains,
the reader can consult  \cite{Kratt1}, \cite{Kratt2} and \cite{IW5}
since, in this paper, we sometimes omit the details and give sketches of proofs.

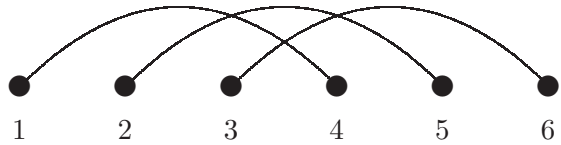
\begin{figure}[b]
\begin{center}
\begin{picture}(200,60)
\put(  0,20){\circle*{8}}
\put( 40,20){\circle*{8}}
\put( 80,20){\circle*{8}}
\put(120,20){\circle*{8}}
\put(160,20){\circle*{8}}
\put(200,20){\circle*{8}}
\put(  0,0){\makebox(0,0)[b]{$1$}}
\put( 40,0){\makebox(0,0)[b]{$2$}}
\put( 80,0){\makebox(0,0)[b]{$3$}}
\put(120,0){\makebox(0,0)[b]{$4$}}
\put(160,0){\makebox(0,0)[b]{$5$}}
\put(200,0){\makebox(0,0)[b]{$6$}}
\qbezier(  0,20)( 60, 80)(120,20)
\qbezier( 40,20)(100, 80)(160,20)
\qbezier( 80,20)(140, 80)(200,20)
\end{picture}
\caption{A perfect matching}\label{figure:matching}
\end{center}
\end{figure}

Let $n$ be a non-negative integer
and assume we are given a $2n$ by $2n$ skew-symmetric matrix $A=(a_{ij})_{1\le i,j\le 2n}$, (i.e. $a_{ji}=-a_{ij}$),
whose entries $a_{ij}$ are in a commutative ring.
The \defterm{Pfaffian} of $A$ is,
by definition,
\begin{equation*}
\label{def_pfaffian}
\Pf(A)=\sum \epsilon(\sigma_{1},\sigma_{2},\hdots,\sigma_{2n-1},\sigma_{2n})\,
a_{\sigma_{1}\sigma_{2}} \dots a_{\sigma_{2n-1}\sigma_{2n}}.
\end{equation*}
where the summation is over all partitions $\{\{\sigma_{1},\sigma_{2}\}_{<},\hdots,\{\sigma_{2n-1},\sigma_{2n}\}_{<}\}$
of $[2n]$ into $2$-elements blocks,
and where $\epsilon(\sigma_{1},\sigma_{2},\hdots,\sigma_{2n-1},\sigma_{2n})$ denotes the sign of the permutation
\begin{equation*}
\begin{pmatrix}
1&2&\cdots&2n\\
\sigma_{1}&\sigma_{2}&\cdots&\sigma_{2n}
\end{pmatrix}.
\end{equation*}
We call a partition $\sigma=\{\{\sigma_{1},\sigma_{2}\}_{<},\hdots,\{\sigma_{2n-1},\sigma_{2n}\}_{<}\}$
of $[2n]$ into $2$-elements blocks \defterm{a perfect matching} or \defterm{$1$-factor} of $[2n]$,
and let ${\cal F}_{n}$ denote the set of all perfect matchings of $[2n]$.
We represent a perfect matching $\sigma$ graphically by
embedding the points $i\in[2n]$ along the $x$-axis in the coordinate plane
and representing each block $\{\sigma_{2i-1},\sigma_{2i}\}_{<}$ by the curve connecting $\sigma_{2i-1}$ to $\sigma_{2i}$
in the upper half plane.
For instance,
the graphical representation of $\sigma=\{\{1,4\},\{2,5\},\{3,6\}\}$ is
the Figure~\ref{figure:matching} bellow.
%%%%%%%%%%%%%%%%%%%%%
% Perfect matching  %
%%%%%%%%%%%%%%%%%%%%%
If we write $\wt(\sigma)=\epsilon(\sigma)\prod_{i=1}^{n}a_{\sigma_{2i-1}\sigma_{2i}}$ for each perfect matching $\sigma$,
then we can restate our definition as
\begin{equation}
\label{eq:Pf_matching}
\Pf(A)=\sum_{\sigma\in{\cal F}_{n}}\wt(\sigma).
\end{equation}
A skew-symmetric matrix $A=(a_{ij})_{1\leq i,j\leq n}$ is uniquely determined by its upper triangular entries
$(a_{ij})_{1\leq i<j\leq n}$.
So we sometimes define a skew-symmetric matrix by describing its upper triangular entries.

Let $O_{m,n}$ denote the $m\times n$ zero matrix and let $E_{n}$ denote the identity matrix $(\delta_{ij})_{1\leq i<j\leq n}$ of size $n$.
Here $\delta_{ij}$ denotes the Kronecker delta.
We use the abbreviation $O_{n}$ for $O_{n,n}$.

For any finite set $S$ and any nonnegative integer $r$,
let $\binom{S}{r}$ denote the set of all $r$-element subsets of $S$.
For example, $\binom{[n]}{r}$ stands for the set of all multi-indices $\{i_1,\dots,i_r\}$
such that $1\leq i_1<\dots<i_r\leq n$.
Let $m$, $n$ and $r$ be integers such that $r\leq m,n$ and let $T$ be an $m$ by $n$ matrix.
For any index sets $I=\{i_1,\dots,i_r \}\in\binom{[m]}{r}$ and $J=\{j_1,\dots,j_r\}\in\binom{[n]}{r}$,
let $\Delta^{I}_{J}(A)$ denote the submatrix obtained by selecting the rows indexed by $I$ and the columns indexed by $J$.
If $r=m$ and $I=[m]$,
we simply write $\Delta_{J}(A)$ for $\Delta^{[m]}_{J}(A)$.
Similarly,
if $r=n$ and $J=[n]$,
we write $\Delta^{I}(A)$ for $\Delta_{[n]}^{I}(A)$.
It is essential that the weight $\omega(\lambda)$ can be expressed by a Pfaffain,
which is a fact proved in \cite{I}:

%----------------------------------------------------------
% Lemma
%----------------------------------------------------------
\begin{theorem}
\label{thm:coefficients}
Let $n$ be a non-negative integer.
Let $\lambda=(\lambda_1,\dots,\lambda_{2n})$ be a partition such that $\ell(\lambda)\leq2n$,
and put $l=(l_1,\dots,l_{2n})=\lambda+\delta_{2n}$.
Define a skew-symmetric matrix $A=(\alpha_{ij})_{i,j\geq0}$
by
\begin{equation*}
\alpha_{ij}=a^{\lceil(j-1)/2\rceil}b^{\lfloor(j-1)/2\rfloor}c^{\lceil i/2\rceil}d^{\lfloor i/2\rfloor}
\end{equation*}
for $i<j$.
Then we have
\begin{equation*}
\Pf\left[\Delta^{I(\lambda)}_{I(\lambda)}\left(A\right)\right]_{1\leq i,j\leq 2n}
=(abcd)^{\binom{n}2}\omega(\lambda),
\end{equation*}
where $I(\lambda)=\{l_{2n},\dots,l_{1}\}$.
\end{theorem}
A variation of this theorem for strict partitions is as follows.
%----------------------------------------------------------
% Lemma (subPfaffians2)
%----------------------------------------------------------
\begin{theorem}
\label{thm:coef_strict}
Let $n$ be a nonnegative integer.
Let $\mu=(\mu_1,\dots,\mu_{n})$ be a strict partition such that
$\mu_1>\dots>\mu_{n}\geq0$.
Let $K(\mu)=\{\mu_{n},\dots,\mu_{1}\}$.
Define a skew-symmetric matrix $B=(\beta_{ij})_{i,j\geq-1}$ by
\begin{equation}
\beta_{ij}
=\begin{cases}
1,
&\text{ if $i=-1$ and $j=0$,}\\
a^{\lceil j/2\rceil}b^{\lfloor j/2\rfloor}z,
&\text{ if $i=-1$ and $j\geq 1$,}\\
%a^{\lceil j/2\rceil}b^{\lfloor j/2\rfloor}
%c^{\lceil i/2\rceil}d^{\lfloor i/2\rfloor}z,
a^{\lceil j/2\rceil}b^{\lfloor j/2\rfloor}z
&\text{ if $i=0$,}\\
a^{\lceil j/2\rceil}b^{\lfloor j/2\rfloor}
c^{\lceil i/2\rceil}d^{\lfloor i/2\rfloor}z^2,
&\text{ if $i>0$,}
\end{cases}
\label{eq:strict}
\end{equation}
for $-1\leq i<j$.
\begin{enumerate}
\item[(i)]
If $n$ is even,
then we have
\begin{equation}
\Pf\left[\Delta^{K(\mu)}_{K(\mu)}\left(B\right)\right]
=\omega(\mu)z^{\ell(\mu)}.
\end{equation}
\item[(ii)]
If $n$ is odd,
then we have
\begin{equation}
\Pf\left[\Delta^{\{-1\}\uplus K(\mu)}_{\{-1\}\uplus K(\mu)}\left(B\right)\right]
=\omega(\mu)z^{\ell(\mu)}.\ \Box
\end{equation}
\end{enumerate}
\end{theorem}

These theorems are easy consequences of the following Lemma which has been proved in \cite[Section~4, Lemma~7]{IW1}.
%----------------------------------------------------------
% Lemma 7 from [IW1]
%----------------------------------------------------------
\begin{lemma}
\label{lemma:product}
Let $x_i$ and $y_j$ be indeterminates,
and let $n$ is a non-negative integer.
Then
\begin{equation*}
\label{prod_pf}
\Pf[x_iy_j]_{1\leq i<j\leq 2n}=\prod_{i=1}^{n}x_{2i-1}\prod_{i=1}^{n}y_{2i}.
\ \Box
\end{equation*}
\end{lemma}

%%%%%%%%%%%%%%%%%%%%%%%%%%%%%%%%%%%%%%%%%%%%%%%%%%%%%%%%%%%%%%%%%%%%%%%%%%%%%
%
% Section 3: Strict Partitions
%
%%%%%%%%%%%%%%%%%%%%%%%%%%%%%%%%%%%%%%%%%%%%%%%%%%%%%%%%%%%%%%%%%%%%%%%%%%%%%

\section{Strict Partitions}\label{Strict_Partitions}

A partition $\mu$ is \defterm{strict} if all its parts are distinct.
One represents the associated shifted diagram of $\mu$
as a diagram in which the $i$th row from the top
has been shifted to the right by $i$ places
so that the first column becomes a diagonal.
A strict partition can be written uniquely
in the form $\mu=(\mu_1,\dots,\mu_{2n})$
where $n$ is an non-negative integer and $\mu_1>\mu_2>\cdots>\mu_{2n}\geq0$.
The \defterm{length} $\ell(\mu)$ is,
by definition,
the number of nonzero parts of $\mu$.
We define the weight function $\omega(\mu)$ exactly the same as in \thetag{\ref{eq:weight}}.
For example,
if $\mu=(8,5,3)$,
then $\ell(\mu)=3$, $\omega(\mu)=a^6b^5c^3d^2$
and its shifted diagram is as follows.
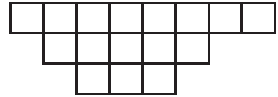
\[
\smyoung{
&&&&&&&\\
\blank&&&&&\\
\blank&\blank&&&\\
}
\]
Let
\begin{equation}
\Psi_N=\Psi_N(a,b,c,d;z)=\sum \omega(\mu)z^{\ell(\mu)},
\label{eq:gen_strict}
\end{equation}
where the sum is over all
strict partitions $\mu$ such that each part of $\mu$ is less than or equal to $N$.
For example,
we have
\begin{align*}
\Psi_0&=1,\\
\Psi_1&=1+az,\\
\Psi_2&=1+a(1+b)z+abcz^2,\\
\Psi_3&=1+a(1+b+ab)z+abc(1+a+ad)z^2+{a}^{3}bcd{z}^{3}.
\end{align*}
In fact, the only strict partition such that $\ell(\mu)=0$ is
$\emptyset$,
the strict partitions $\mu$ such that $\ell(\mu)=1$
and $\mu_1\leq3$ are the following three:
\[
\smyoung{
a\\
\blank\\}
\qquad
\smyoung{
a&b\\
\blank\\}
\qquad
\smyoung{
a&b&a\\
\blank\\}\ ,
\]
the strict partitions $\mu$ such that $\ell(\mu)=2$
and $\mu_1\leq3$ are the following three:
\[
\smyoung{
a&b\\
\blank&c\\}
\qquad
\smyoung{
a&b&a\\
\blank&c\\}
\qquad
\smyoung{
a&b&a\\
\blank&c&d\\}\ ,
\]
and the strict partition $\mu$ such that $\ell(\mu)=3$ and $\mu_1\leq3$ is
the following one:
\[
\smyoung{
a&b&a\\
\blank&c&d\\
\blank&\blank&a\\
}\ .
\]
The sum of the weights of these strict partitions is equal to $\Psi_3$.
In this section
we always assume $|a|,|b|,|c|,|d|<1$.
One of the main results of this section is that
the even terms and the odd terms of $\Psi_N$ respectively satisfy the associated Al-Salam-Chihara recurrence relation:
%----------------------------------------------------------
% Main Theorem
%----------------------------------------------------------
\begin{theorem}
\label{th:aac}
Set $q=abcd$.
Let $\Psi_{N}=\Psi_{N}(a,b,c,d;z)$ be as in \thetag{\ref{eq:gen_strict}}
and put $X_N=\Psi_{2N}$ and $Y_N=\Psi_{2N+1}$.
Then $X_N$ and $Y_N$ satisfy
\begin{align}
X_{N+1}&= \left\{1+ab+a(1+bc)z^2q^N\right\}X_N
\nonumber\\
&\qquad\qquad\qquad-ab(1-z^2q^N)(1-acz^2q^{N-1})X_{N-1},
\label{eq:recX}\\
Y_{N+1}&= \left\{1+ab+abc(1+ad)z^2q^N\right\}Y_N
\nonumber\\
&\qquad\qquad\qquad-ab(1-z^2q^N)(1-acz^2q^{N})Y_{N-1},
\label{eq:recY}
\end{align}
where
$X_0=1$, $Y_0=1+az$, $X_1=1+a(1+b)z+abcz^2$ and
\[
Y_1=1+a(1+b+ab)z+abc(1+a+ad)z^2+a^{3}bcdz^3.
\]
Especially,
if we put $X_N'=(ab)^{-\frac{N}2}X_{N}$ and $Y_N'=(ab)^{-\frac{N}2}Y_{N}$,
then $X_N'$ and $Y_N'$ satisfy
\begin{align}
\left\{(ab)^{\frac12}+(ab)^{-\frac12}\right\}X_{N}'
&= X_{N+1}'
-a^{\frac12}b^{-\frac12}(1+bc)z^2q^N X_N'
\nonumber\\
&\qquad+(1-z^2q^N)(1-acz^2q^{N-1})X_{N-1}',
\label{eq:recXX}\\
\left\{(ab)^{\frac12}+(ab)^{-\frac12}\right\}Y_{N}'
&= Y_{N+1}'
-a^{\frac12}b^{\frac12}c(1+ad)z^2q^N Y_N'
\nonumber\\
&\qquad+(1-z^2q^N)(1-a^2bc^2dz^2q^{N-1})Y_{N-1}',
\label{eq:recYY}
\end{align}
where
$X_0'=1$, $Y_0'=1+az$, $X_1'=(ab)^{-\frac12}+a^{\frac12}b^{-\frac12}(1+b)z+(ab)^{\frac12}cz^2$ and
\[
Y_1'=(ab)^{-\frac12}+a^{\frac12}b^{-\frac12}(1+b+ab)z+a^{\frac12}b^{\frac12}c(1+a+ad)z^2+a^{\frac52}b^{\frac12}cdz^3.
\]
\end{theorem}

Thus \thetag{\ref{eq:recXX}} agrees with the associated Al-Salam-Chihara recurrence relation \thetag{\ref{eq:AASC}}
where $u=a^{\frac12}b^{\frac12}$,
$\alpha=-a^{\frac12}b^{\frac12}c$,
$\beta=-a^{\frac12}b^{-\frac12}$
 and
$t=z^2$,
and \thetag{\ref{eq:recYY}} also agrees with \thetag{\ref{eq:AASC}}
where $u=a^{\frac12}b^{\frac12}$,
$\alpha=-a^{\frac12}b^{\frac12}c$,
$\beta=-a^{\frac32}b^{\frac12}cd$
 and $t=z^2$.
One concludes that,
when $|a|,|b|,|c|,|d|<1$,
the solutions of \thetag{\ref{eq:recX}} and \thetag{\ref{eq:recY}} are expressed by the linear combinations of \thetag{\ref{eq:Q1}} and \thetag{\ref{eq:Q2}} as follows.
%----------------------------------------------------------
% Main Theorem
%----------------------------------------------------------
\begin{theorem}
\label{th:gen_strict}
Assume $|a|,|b|,|c|,|d|<1$ and
set $q=abcd$.
Let $\Psi_{N}=\Psi_{N}(a,b,c,d;z)$ be as in \thetag{\ref{eq:gen_strict}}.
\begin{enumerate}
%-------------------------------------------%
% General X_n
%-------------------------------------------%
\item[(i)]
Put $X_N=\Psi_{2N}$.
Then we have
\begin{align}
X_{N}
&=
\frac
{(-a z^2 q,-abc;q)_{\infty}}
{(-a,-abcz^2;q)_{\infty}}
\nonumber\\
&\qquad\times\left\{
(s^{X}_{0}X_{1} - s^{X}_{1}X_{0})
(-abc z^2;q)_N\,
\PHI{q^{-N}z^{-2}}{-b^{-1}}{-(abc)^{-1}q^{-N+1}z^{-2}}{-c^{-1}q}
\right.\nonumber\\
&\qquad\left.
+(r^{X}_{1}X_{0} - r^{X}_{0}X_{1})
(ab)^{N}\,
\frac{(qz^2,ac z^2 ;q)_N}{(-a q z^2;q)_N}\,
\PHI{q^{N+1} z^2}{-c^{-1} q}{-a q^{N+1} z^2}{-abc}
\right\},
\label{eq:solX}
\end{align}
where
\begin{align*}
r^{X}_{0}
&=
\PHI{z^{-2}}{-b^{-1}}{-(abc)^{-1}z^{-2}q}{-c^{-1}q},
\\
s^{X}_{0}
&=
\PHI{ z^2 q}{-c^{-1}q}{-a z^2 q}{-abc},
\\
r^{X}_{1}
&=(1+abc z^2)\,
\PHI{z^{-2}q^{-1}}{-b^{-1}}{-(abc)^{-1}z^{-2}}{-c^{-1}q},
\\
s^{X}_{1}
&=
\frac{ab(1-z^2q)(1-ac z^2)}{1+a z^2 q}\,
\PHI{ z^2 q^{2}}{-c^{-1}q}{-a z^2 q^{2}}{-abc}.
\end{align*}
%-------------------------------------------%
% General Y_n
%-------------------------------------------%
\item[(ii)]
Put $Y_N=\Psi_{2N+1}$.
Then we have
\begin{align}
Y_{N}
&=
\frac
{(-a q^2 z^2,-abc;q)_{\infty}}
{(-a q,-abc z^2;q)_{\infty}}
\nonumber\\
&\times\biggl\{(s^{Y}_{0}Y_{1} - s^{Y}_{1}Y_{0})
(-abc z^2;q)_N\,
\PHI{q^{-N}z^{-2}}{-acd}{-(abc)^{-1}q^{-N+1}z^{-2}}{-c^{-1}q}
\nonumber\\
&+(r^{Y}_{1}Y_{0} - r^{Y}_{0}Y_{1})
(ab)^{N}\,
\frac{(q z^2,acq z^2 ;q)_N}{(-a q^2 z^2;q)_N}\,
\PHI{q^{N+1} z^2}{-c^{-1}q}{-a q^{N+2} z^2}{-abc}
\biggr\},
\label{eq:solY}
\end{align}
where
\begin{align*}
r^{Y}_{0}
&=
\PHI{z^{-2}}{-acd}{(-abc)^{-1}qz^{-2}}{-c^{-1}q},
\\
r^{Y}_{1}
&=(1+abc z^2)\,
\PHI{q^{-1}z^{-2}}{-acd}{-(abc)^{-1}z^{-2}}{-c^{-1}q},
\\
s^{Y}_{0}
&=
\PHI{z^2q}{-c^{-1}q}{-a q^2 z^2}{-abc},
\\
s^{Y}_{1}
&=
\frac{ab(1-z^2q)(1-acq z^2)}{1+a q^2 z^2}\,
\PHI{ z^2 q^{2}}{-c^{-1}q}{-a q^{3} z^2}{-abc}.
\end{align*}
\end{enumerate}
\end{theorem}
%----------------------------------------------------------
% Corollary (Limit)
%----------------------------------------------------------
If we take the limit $N\to\infty$ in \thetag{\ref{eq:solX}} and \thetag{\ref{eq:solY}},
then by using \thetag{\ref{eq:limit}},
we obtain the following generalization of Boulet's result
(see Corollary~\ref{boulet:strict}).
\begin{corollary}
\label{cor:limit_strict}
Assume $|a|,|b|,|c|,|d|<1$ and
set $q=abcd$.
Let $s^{X}_{i}$, $s^{Y}_{i}$, $X_{i}$, $Y_{i}$ ($i=0,1$) be as in the above theorem.
Then we have
\begin{align}
\sum_{\mu}\omega(\mu)z^{\ell(\mu)}
&=
\frac{(-abc , -a z^2 q;q)_{\infty}}
{(ab;q)_{\infty}}(s^{X}_{0}X_{1} - s^{X}_{1}X_{0})
\nonumber\\
&=
\frac{( -abc , -a z^2 q^2;q)_{\infty}}
{(ab;q)_{\infty}}(s^{Y}_{0}Y_{1} - s^{Y}_{1}Y_{0}),
\label{eq:limit_strict}
\end{align}
where the sum runs over all strict partitions and the first terms
 are as follows:
\begin{align*}
1+\frac{a(1+b)}{1-ab}z
+\frac{abc(1+a+ad+abd)}{(1-ab)(1-q)}z^2
+\frac{a^2q(1+b)(1+bc+abc+bq)}{(1-ab)(1-q)(1-abq)}z^3
+O(z^4).
\end{align*}
\end{corollary}
%----------------------------------------------------------
% Main Corollary
%----------------------------------------------------------
On the other hand, by plugging
 $z=1$ into \thetag{\ref{eq:solX}} and \thetag{\ref{eq:solY}},
we conclude that the solutions of the recurrence relations \thetag{\ref{eq:recXX}} and \thetag{\ref{eq:recYY}}
with the above initial condition
are exactly the Al-Salam-Chihara polynomials, which give two finite versions of Boulet's result.
\begin{corollary}
\label{hyper:strict}
Put $u=\sqrt{ab}$, $x=\frac{u+u^{-1}}2$ and $q=abcd$.
Let $\Psi_{N}(a,b,c,d;z)$ be as in \thetag{\ref{eq:gen_strict}}.
\begin{enumerate}
\item[(i)]
The polynomial $\Psi_{2N}(a,b,c,d;1)$ is given by
\begin{align}
\Psi_{2N}(a,b,c,d;1)
&=
(ab)^{\frac{N}2}Q_{N}(x;-a^{\frac12}b^{\frac12}c,-a^{\frac12}b^{-\frac12}|q),
\nonumber\\
&=(-a;q)_{N}\,
\PHI{q^{-N}}{-c}{-a^{-1}q^{-N+1}}{-bq}.
\label{eq:st_al1}
\end{align}
\item[(ii)]
The polynomial $\Psi_{2N+1}(a,b,c,d;1)$ is given by
\begin{align}
\Psi_{2N+1}(a,b,c,d;1)
&=
(1+a)(ab)^{\frac{N}2}
Q_{N}(x;-a^{\frac12}b^{\frac12}c,-a^{\frac32}b^{\frac12}cd|q)
\nonumber\\
&=(-a;q)_{N+1}\,
\PHI{q^{-N}}{-c}{-a^{-1}q^{-N}}{-b}.
\label{eq:st_al2}
\end{align}
\end{enumerate}
\end{corollary}
Substituting $a=zyq$, $b=z^{-1}yq$, $c=zy^{-1}q$ and $d=z^{-1}y^{-1}q$
into Theorem~\ref{hyper:strict}(see \cite{Bou}),
then we immediately obtain the strict version of Andrews' result
(see Corollary~\ref{cor:andrews}).
%----------------------------------------------------------
% Corollary
%----------------------------------------------------------
\begin{corollary}
\label{cor:strict_andrews}
\begin{equation}
\sum_{{\mu\text{ strict partitions}}\atop{\mu_1\leq 2N}}
z^{{\cal O}(\mu)}y^{{\cal O}(\mu')}q^{|\mu|}
=\sum_{j=0}^{N}
\left[{{N}\atop{j}}\right]_{q^4}
(-zyq;q^4)_{j}(-zy^{-1}q;q^4)_{N-j}(yq)^{2N-2j},
\label{eq:strict_andrews01}
\end{equation}
and
\begin{equation}
\sum_{{\mu\text{ strict partitions}}\atop{\mu_1\leq 2N+1}}
z^{{\cal O}(\mu)}y^{{\cal O}(\mu')}q^{|\mu|}
=\sum_{j=0}^{N}
\left[{{N}\atop{j}}\right]_{q^4}
(-zyq;q^4)_{j+1}(-zy^{-1}q;q^4)_{N-j}(yq)^{2N-2j},
\label{eq:strict_andrews02}
\end{equation}
where
\[
\left[{{N}\atop{j}}\right]_{q}
=\begin{cases}
\frac{(1-q^{N})(1-q^{N-1})\cdots(1-q^{N-j+1})}
{(1-q^{j})(1-q^{j-1})\cdots(1-q)},
&\text{ for $0\leq j\leq N$,}
\\
0,
&\text{ if $j<0$ and $j>N$.}
\end{cases}
\]
\end{corollary}
Letting $N\rightarrow\infty$ in Theorem~\ref{hyper:strict}
or  setting $z=1$ in \thetag{\ref{eq:limit_strict}},
 we  obtain the following result of Boulet
(cf. \cite[Corollary~2]{Bou}).
%----------------------------------------------------------
% Theorem (Boulet) strict partitions
%----------------------------------------------------------
\begin{corollary}
\label{boulet:strict}
(Boulet)
Let $q=abcd$,
then
\begin{align}
\sum_{\mu}\omega(\mu)
=
\frac{(-a;q)_{\infty}(-abc;q)_{\infty}}
{(ab;q)_{\infty}},
\label{eq:boulet2}
\end{align}
where the sum runs over all strict partitions.
\end{corollary}

To prove Theorem~\ref{th:aac},
we need several steps.
Our strategy is as follows: write
the weight $\omega(\mu)z^{\ell(\mu)}$ as a Pfaffain
(Theorem~\ref{thm:coef_strict})
and apply the minor summation formula (Lemma~\ref{lem:pf_sum})
to make the sum of the weights into a single Pfaffian
(Theorem~\ref{th:pf_strict}).
Then we make use of the Pfaffian to derive a recurrence relation
(Proposition~\ref{prop:rec_dist}).
We also give another proof of the recurrence relation by a combinatorial argument
(Remark~\ref{rem:comb}).

Let $J_n$ denote the square matrix of size $n$ whose $(i,j)$th entry is $\delta_{i,n+1-j}$.
We simply write $J$ for $J_n$ when there is no fear of confusion on the size $n$.
We need the following result on a sum of Pfaffians~\cite[Theorem of Section~4]{TW}.
\begin{lemma}
\label{lem:pf_sum}
Let $n$ be a positive integer.
Let $A=(a_{ij})_{1\leq i,j\leq n}$ and $B=(b_{ij})_{1\leq i,j\leq n}$
be skew symmetric matrices of size $n$.
Then
\begin{align}
&\sum_{t=0}^{\lfloor n/2\rfloor}z^t
\sum_{I\in\binom{[n]}{2t}}\gamma^{|I|}
\Pf\left(\Delta^I_I(A)\right)
\Pf\left(\Delta^I_I(B)\right)
=\Pf\begin{bmatrix}
J_n\,{}^t\kern-2pt AJ_n&J_n\\
-J_n&C
\end{bmatrix}
%=\Pf\left(\Pf(A)^{-1}\widehat A+C\right)
,
\label{eq:pf_sum}
\end{align}
where $|I|=\sum_{i\in I}i$ and $C=(C_{ij})_{1\leq i,j\leq n}$ is given by
$C_{ij}=\gamma^{i+j}b_{ij}z$.
\end{lemma}
This lemma is a special case of Lemma~\ref{lem:pf_sum2}, so
 a proof will be given later.
\bigbreak

Let $S_{n}$ denote
the $n\times n$ skew-symmetric matrix
whose $(i,j)$th entry is $1$ for $0\leq i<j\leq n$.
As a corollary of Lemma~\ref{lem:pf_sum},
we obtain the following expression of the sum of the weight $\omega(\mu)$  by a single Pfaffian.
%----------------------------------------------------------
% Theorem
%----------------------------------------------------------
\begin{theorem}
\label{th:pf_strict}
Let $N$ be a nonnegative integer.
\begin{align}
\Psi_{N}(a,b,c,d;z)
=\Pf\begin{bmatrix}
S_{N+1}&J_{N+1}\\
-J_{N+1}&B
\end{bmatrix},
\label{eq:pf_dist}
\end{align}
where $B=(\beta_{ij})_{0\leq i<j\leq N}$ is
the $(N+1)\times(N+1)$ skew-symmetric matrix
whose $(i,j)$th entry $\beta_{ij}$ is defined in \thetag{\ref{eq:strict}}.
%***% \[
%***% \beta_{ij}=\begin{cases}
%***% a^{\lceil j/2\rceil}
%***% b^{\lfloor j/2\rfloor}
%***% c^{\lceil i/2\rceil}
%***% d^{\lfloor i/2\rfloor}z
%***% &\text{ if $0=i<j\leq n$,}
%***% \\
%***% a^{\lceil j/2\rceil}
%***% b^{\lfloor j/2\rfloor}
%***% c^{\lceil i/2\rceil}
%***% d^{\lfloor i/2\rfloor}z^2
%***% &\text{ if $0<i<j\leq n$.}
%***% \end{cases}
%***% \]
\end{theorem}
\begin{demo}{Proof}
By Theorem~\ref{thm:coef_strict},
we have
\[
\Psi_{N}(a,b,c,d;z)
=\sum_{t=0}^{\lfloor (N+1)/2\rfloor}
\sum_{{\mu=(\mu_1,\dots,\mu_{2t})}
\atop{N\geq\mu_1>\cdots>\mu_{2t}\geq0}}
\Pf\left(
\Delta^{K(\mu)}_{K(\mu)}\left(B\right)
\right).
\]
If we take $S=(1)_{0\leq i<j\leq N}$,
then we have $\Pf\left(\Delta^I_I(S)\right)=1$ for any subset $I\subseteq[0,N]$ of even cardinality.
(For detailed arguments on sub-pfaffains, see \cite{IW5}).
Thus \thetag{\ref{eq:pf_dist}} follows from Lemma~\ref{lem:pf_sum}.
$\Box$
\end{demo}
For example,
if $N=3$,
then the Pfaffian in the right-hand side of \thetag{\ref{eq:pf_dist}} looks
\begin{equation*}
\Pf
\left[\begin{array}{cccc|cccc}
0&1&1&1&0&0&0&1\\
-1&0&1&1&0&0&1&0\\
-1&-1&0&1&0&1&0&0\\
-1&-1&-1&0&1&0&0&0\\\hline
0&0&0&-1&0&az&abz&{a}^{2}bz\\
0&0&-1&0&-az&0&abcz^{2}&{a}^{2}bcz^{2}\\
0&-1&0&0&-abz&-abcz^{2}&0&{a}^{2}bcdz^{2}\\
-1&0&0&0&-{a}^{2}bz&-{a}^{2}bcz^{2}&-{a}^{2}bcdz^{2}&0
\end{array}\right],
\end{equation*}
and this is equal to
$
1+a(1+b+ab)z+abc(1+a+ad)z^{2}+{a}^{3}bcd{z}^{3}
$.

By performing elementary transformations on rows
and columns of the matrix,
we obtain the following recurrence relation:
%----------------------------------------------------------
% Proposition
%----------------------------------------------------------
\begin{prop}
\label{prop:rec_dist}
Let $\Psi_{N}=\Psi_{N}(a,b,c,d;z)$ be as above.
Then we have
\begin{align}
&\Psi_{2N}=(1+b)\Psi_{2N-1}+(a^{N}b^{N}c^{N}d^{N-1}z^{2}-b)\Psi_{2N-2},
\label{eq:st1}
\\
&\Psi_{2N+1}=(1+a)\Psi_{2N}+(a^{N+1}b^{N}c^{N}d^{N}z^{2}-a)\Psi_{2N-1},
\label{eq:st2}
\end{align}
for any positive integer $N$.
\end{prop}
%----------------------------------------------------------
% Proof of Proposition
%----------------------------------------------------------
\begin{demo}{Proof}
Let $A$ denote the $2(N+1)\times2(N+1)$ skew symmetric matrix
$\begin{bmatrix}
S_{N+1}&J_{N+1}\\
-J_{N+1}&B
\end{bmatrix}$ as on the right-hand side of \thetag{\ref{eq:pf_dist}}.
Here we assume row/column indices start at 0.
So,
for example,
the row indices for the upper $(N+1)$ rows are $i$, $i=0,\dots,N$,
and the row indices for the lower $(N+1)$ rows are $i+N+1$, $i=0,\dots,N$.
%We perform certain elementary transformations on $A$ and change it into a simpler form.
Now, subtract $a$ times $(j+N)$th column from $(j+N+1)$th column if $j$ is odd,
or subtract $b$ times $(j+N)$th column from $(j+N+1)$th column if $j$ is even,
for $j=N,N-1,\dots,1$.
To make our matrix skew-symmetric,
subtract $a$ times $(i+N)$th row from $(i+N+1)$th row if $i$ is odd,
or subtract $b$ times $(i+N)$th row from $(i+N+1)$th row if $i$ is even,
for $i=N,N-1,\dots,1$.
Next subtract $(i+1)$th row from $i$th row for $i=0,1,\dots,N-1$,
then we also subtract $(j+1)$th column from $j$th column for $j=0,1,\dots,N-1$.
By these transformations,
we obtain a skew symmetric matrix $A'=\begin{bmatrix}
P&Q\\
-{}^tQ&R
\end{bmatrix}$,
where $P=(\delta_{i+1,j})_{0\leq i<j\leq N}$,
$Q=(q_{ij})_{0\leq i<j\leq N}$
and
$R=(r_{ij})_{0\leq i<j\leq N}$
are given by
\begin{align*}
&q_{ij}=\begin{cases}
-1
&\text{ if $i+j=N-1$,}\\
1
&\text{ if $i=N$ and $j=0$,}\\
1+a^{\chi(\text{$j$ is odd})}b^{\chi(\text{$j$ is even})}
&\text{ if $i+j=N$ and $j\geq1$,}\\
-a^{\chi(\text{$j$ is odd})}b^{\chi(\text{$j$ is even})}
&\text{ if $i+j=N+1$,}\\
0
&\text{ otherwise,}
\end{cases}
\\
&r_{ij}=\begin{cases}
az\delta_{1,j}
&\text{ if $i=0$,}
\\
a^{\lceil (i+1)/2\rceil}b^{\lfloor (i+1)/2\rfloor}c^{\lceil i/2\rceil}d^{\lfloor i/2\rfloor}z^2\delta_{i+1,j}
&\text{ if $i>0$.}
\end{cases}
\end{align*}
Here $\chi(A)$ stands for $1$ if the statement $A$ is true and $0$ otherwise.
For example,
if $N=3$,
then $A'$ looks as follows:
\[
\left[ \begin {array}{cccc|cccc}
0&1&0&0&0&0&-1&1+a\\
-1&0&1&0&0&-1&1+b&-a\\
0&-1&0&1&-1&1+a&-b&0\\
0&0&-1&0&1&-a&0&0\\\hline
0&0&1&-1&0&az&0&0\\
0&1&-1-a&a&-az&0&abcz^{2}&0\\
1&-1-b&b&0&0&-abcz^{2}&0&{a}^{2}bcdz^{2}\\
-1-a&a&0&0&0&0&-{a}^{2}bcdz^{2}&0
\end {array}\right].
\]
By expanding $\Pf(A')$ along the first row/column,
we obtain the desired formula.
$\Box$
\end{demo}

\begin{remark}
\label{rem:comb}
Proposition~\ref{prop:rec_dist} can be also proved by a combinatorial argument as follows.
\end{remark}

\begin{demo}{Combinatorial proof of Proposition~\ref{prop:rec_dist}}
By definition, the generating function for strict partitions $\mu=(\mu_1,\mu_2,\dots)$
such that
$\mu_1=2N$
and $\mu_{2}\leq 2N-2$ is equal to
$$
b(\Psi_{2N-1}-\Psi_{2N-2}).
$$
That for strict partitions such that $\mu_1=2N$
and $\mu_{2}=2N-1$ is equal to
$$
a^{N}b^{N}c^{N}d^{N-1}z^{2}\Psi_{2N-2}.
$$
Finally the generating function of strict partitions such that $\mu_1\leq 2N-1$ is equal to
$\Psi_{2N-1}$.
Summing up we get \thetag{\ref{eq:st1}}.
The same argument works to prove \thetag{\ref{eq:st2}}.
$\Box$
\end{demo}

Note that
one can immediately derive Theorem~\ref{th:aac} from Proposition~\ref{prop:rec_dist}
by substitution.
Thus, if one use \thetag{\ref{eq:Q_n}},
then he immediately derive Theorem~\ref{th:gen_strict}
by a simple computation.
%----------------------------------------------------------
% Proof of Theorem 3.2
%----------------------------------------------------------
\begin{demo}{Proof of Theorem~\ref{th:gen_strict}}
%***% \begin{demo}{Proof of Theorem~\ref{hyper:strict}}
%***% Let $X_{n}=(ab)^{\frac{n}2}X_{n}'$.
%***% Then \thetag{\ref{eq:recX}} is rewritten as
%***% \begin{equation}
%***% \{(ab)^{\frac12}+(ab)^{-\frac12}\}X_{n}'=X_{n+1}'
%***% -\{a^{\frac12}b^{-\frac12}(1+bc)z^2q^{n}\}X_{n}'
%***% +(1-z^2q^n)(1-acz^2q^{n-1})X_{n-1}'
%***% \label{eq:rec_alX}
%***% \end{equation}
%***% Similarly,
%***% if we put $Y_{n}=(ab)^{\frac{n}2}Y_{n}'$,
%***% then, by \thetag{\ref{eq:recY}}, we obtain
%***% \begin{equation}
%***% \{(ab)^{\frac12}+(ab)^{-\frac12}\}Y_{n}'=Y_{n+1}'
%***% -\{a^{\frac12}b^{\frac12}c(1+ad)z^2q^{n}\}Y_{n}'
%***% +(1-z^2q^n)(1-a^2bc^2dz^2q^{n-1})Y_{n-1}'
%***% \label{eq:rec_alY}
%***% \end{equation}
%***% When $z=1$,
%***%  \thetag{\ref{eq:recX} with
%***% \thetag{\ref{eq:AASC}}, then we conclude that this is the
%***% associated Al-Salam-Chihara recurrence relation with
%***% $u=a^{\frac12}b^{\frac12}$, $\alpha=-a^{\frac12}b^{-\frac12}$,
%***% $\beta=-a^{\frac12}b^{\frac12}c$ and $t=z^2$.
%***% Hereafter we assume
%***% $z=1$.
%***% Assume $z=1$ hereafter.
%***% Computing $X_{-1}$ using the relation backward with
%***% $X_{0}=1$ and $X_{1}=1+a(1+b+bc)$, we obtain $X_{-1}=0$. Thus,
%***% $X_{N}$ is expressed by ${}_2\phi_1$, and an easy computation
%***% leads to \thetag{\ref{eq:st_al1}}.
Let $u=\sqrt{ab}$, $t=z^2$ and $q=abcd$.
By \thetag{\ref{eq:recXX}},
$X_{N}'$ satisfies
the associated Al-Salam-Chihara recurrence relation
\thetag{\ref{eq:AASC}}
with
$\alpha=-a^{\frac12}b^{\frac12}c$ and $\beta=-a^{\frac12}b^{-\frac12}$.
Note that $|u|<1$ and $|q|<|\alpha|<1$ hold.
Thus,
by \thetag{\ref{eq:Q_n}},
we conclude that $X_N$ is given by \thetag{\ref{eq:solX}}.
%***% By Theorem~\ref{th:aac},
A similar argument shows that
$Y_{N}'$ satisfies \thetag{\ref{eq:AASC}} with
$\alpha=-a^{\frac32}b^{\frac12}c$ and
$\beta=-a^{\frac12}b^{\frac12}cd$,
which implies $Y_N$ is given by \thetag{\ref{eq:solY}}.
%***% Computing $Y_{-1}$ using the relation backward again,
%***% we have $Y_{-1}=0$ and $Y_{0}=1+a$. Thus we obtain
%***% \thetag{\ref{eq:st_al2}}.
$\Box$
\end{demo}

%----------------------------------------------------------
% Proof of Theorem 3.4
%----------------------------------------------------------
\begin{demo}{Proof of Theorem~\ref{hyper:strict}}
First,
if we set $z=1$ in \thetag{\ref{eq:solX}},
we obtain that
 $s^{X}_{0}X_{1} - s^{X}_{1}X_{0}=1+a$ and $r^{X}_{1}X_{0} - r^{X}_{0}X_{1}=0$
which immediately imply
\begin{equation*}
X_{N}=
(-abc;q)_{N}\,
\PHI{q^{-N}}{-b^{-1}}{-(abc)^{-1}q^{-N+1}}{-c^{-1}q}.
\end{equation*}
In fact,
it is easy to see that,
when $z=1$,
$r^{X}_{0}=0$ and $r^{X}_{1}=1+abc+a(1+b)$,
which immediately implies $r^{X}_{1}X_{0} - r^{X}_{0}X_{1}=0$.
A similar computation shows that,
if we put $z=1$,
we have
\begin{align*}
s^{X}_{0}
&=
\sum_{n=0}^{\infty}
\frac{(-c^{-1}q;q)_{n}}{(-aq;q)_{n}}
(-abc)^n,
\\
s^{X}_{1}
&=
ab(1-ac)
\sum_{n=0}^{\infty}
\frac{(1-q^{n+1})(-c^{-1}q;q)_{n}}{(-aq;q)_{n+1}}
(-abc)^n.
\end{align*}
If we use these equalities,
then we obtain
\begin{align*}
s^{X}_{0}X_{1} - s^{X}_{1}X_{0}
&=
(1+a)\sum_{n=0}^{\infty}
\frac{(-c^{-1}q;q)_{n}}{(-aq;q)_{n+1}}(-abc)^{n}
\{a+abc+a(1+b)q^{n+1}\}
\\
&=
(1+a)\left\{
\sum_{n=0}^{\infty}
\frac{(-c^{-1}q;q)_{n}}{(-aq;q)_{n}}(-abc)^{n}
-
\sum_{n=0}^{\infty}
\frac{(-c^{-1}q;q)_{n+1}}{(-aq;q)_{n+1}}(-abc)^{n+1}
\right\}.
\end{align*}
Thus the right-hand side equals $1$,
and this proves \thetag{\ref{eq:st_al1}}.
By a similar argument
we can derive \thetag{\ref{eq:st_al2}} from \thetag{\ref{eq:solY}}.
The details are left to the reader.
$\Box$
\end{demo}

%----------------------------------------------------------
% Proof of Corollary
%----------------------------------------------------------
\begin{demo}{Proof of Corollary~\ref{cor:strict_andrews}}
We first claim that
\begin{equation}
\label{eq:key_andrews}
\Psi_{2N}(a,b,c,d;1)=\sum_{k=0}^{N}
\left[{{N}\atop{k}}\right]_q
(-a;q)_{k}(-c;q)_{N-k}(ab)^{N-k}.
\end{equation}
Then \thetag{\ref{eq:strict_andrews01}} is an easy consequence of
\thetag{\ref{eq:key_andrews}} by substituting $a\leftarrow zyq$, $b\leftarrow z^{-1}yq$, $c\leftarrow zy^{-1}q$ and $d\leftarrow z^{-1}y^{-1}q$.
In fact,
using $(q^{-N};q)_{k}=\frac{(q;q)_{N}}{(q;q)_{N-k}}(-1)^k q^{\binom{k}2-Nk}$,
we have
\begin{align*}
\PHI{q^{-N}}{-c}{-a^{-1}q^{-N+1}}{-bq}
=\sum_{k=0}^{N}
\left[{{N}\atop{k}}\right]_q
\frac{(-c;q)_{N-k}}{(-a^{-1}q^{-N+1};q)_{N-k}}
q^{\binom{N-k}2-N(N-k)}(bq)^{N-k}.
\end{align*}
Substitute $(-a^{-1}q^{-N+1};q)_{N-k}=\frac{(-a;q)_{N}}{(-a;q)_{k}}a^{-N+k}q^{-\binom{N}2+\binom{k}2}$
into this identity to show that the right-hand side equals
\begin{align*}
\sum_{k=0}^{N}
\left[{{N}\atop{k}}\right]_q
\frac{(-a;q)_{k}(-c;q)_{N-k}}{(-a;q)_{N}}
(ab)^{N-k}.
\end{align*}
Finally, use \thetag{\ref{eq:st_al1}} to obtain \thetag{\ref{eq:key_andrews}}.
The proof of \thetag{\ref{eq:strict_andrews02}} reduces to
\begin{equation}
\label{eq:key_andrews2}
\Psi_{2N+1}(a,b,c,d;1)=\sum_{k=0}^{N}
\left[{{N}\atop{k}}\right]_q
(-a;q)_{k+1}(-c;q)_{N-k}(ab)^{N-k},
\end{equation}
which is derived from \thetag{\ref{eq:st_al2}} similarly.
$\Box$
\end{demo}

%----------------------------------------------------------
% Proof of Corollary (Boulet)
%----------------------------------------------------------
\begin{demo}{Proof of Corollary~\ref{boulet:strict}}
%By \thetag{\ref{eq:key_andrews}}, we have
%\begin{align*}
%\Psi_{2N}(a,b,c,d;1)
%&=\sum_{k=0}^{N}
%\left[{{N}\atop{k}}\right]_q
%(-a;q)_{N-k}(-c;q)_{k}(ab)^{k}\\
%&\rightarrow(-a;q)_{\infty}\sum_{k=0}^{\infty}
%\frac{(-c;q)_{k}}{(q;q)_{k}}(ab)^{k}
%\end{align*}
%when $N\rightarrow\infty$.
By replacing $k$ by $N-k$ and letting $N$ to $+\infty$,
we get
\begin{align*}
\lim_{N\rightarrow\infty}\Psi_{2N}(a,b,c,d;1)
=(-a;q)_{\infty}\sum_{k=0}^{\infty}
\frac{(-c;q)_{k}}{(q;q)_{k}}(ab)^{k}
=\frac{(-a;q)_{\infty}(-abc;q)_{\infty}}{(ab;q)_{\infty}}
\end{align*}
%Use the $q$-binomial theorem
%$\sum_{k=0}^{\infty}\frac{(a;q)_{k}}{(q;q)_{k}}z^{k}=\frac{(az;q)_{\infty}}{(z;q)_{\infty}}$
%to obtain $\lim_{N\rightarrow\infty}\Psi_{2N}(a,b,c,d;1)=\frac{(-a;q)_{\infty}(-abc;q)_{\infty}}{(ab;q)_{\infty}}$.
where the last equality follows from the $q$-binomial formula (see \cite{GR}).
Similarly
we can derive the limit
from \thetag{\ref{eq:key_andrews2}}.

Note that we can also derive \thetag{\ref{eq:boulet2}}
from \thetag{\ref{eq:limit_strict}}
by the same argument as in the proof of Theorem~\ref{hyper:strict}.
$\Box$
\end{demo}

%%%%%%%%%%%%%%%%%%%%%%%%%%%%%%%%%%%%%%%%%%%%%%%%%%%%%%%%%%%%%%%%%%%
%
% Section 4: Ordinary Partitions
%
%%%%%%%%%%%%%%%%%%%%%%%%%%%%%%%%%%%%%%%%%%%%%%%%%%%%%%%%%%%%%%%%%%%

\section{Ordinary Partitions}

First we present a generalization of Andrews' result in \cite{And}.
Let us consider
\begin{equation}
\Phi_{N}
=\Phi_{N}(a,b,c,d;z)
=\sum_{{\lambda}\atop{\lambda_1\leq N}}\omega(\lambda)z^{\ell(\lambda)},
\label{eq:gen_ord}
\end{equation}
where the sum runs over all partitions $\lambda$
such that each part of $\lambda$ is less than or equal to $N$.
For example,
the first few terms can be computed directly as follows:
\begin{align*}
&\Phi_{0}=1,\\
&\Phi_{1}=\frac{1+az}{1-acz^2},\\
&\Phi_{2}=\frac{1+a(1+b)z+abcz^2}{(1-acz^2)(1-qz^2)},\\
&\Phi_{3}=\frac{1+a(1+b+ab)z+abc(1+a+ad)z^2+a^3bcdz^3}
{(1-z^2ac)(1-z^2q)(1-z^2acq)},
%\\
%&\Phi_{4}=\frac{1+za(1+b)(1+ab)+z^2a(bc+abc+ab^2c+q+bq+bcq)+z^3a^2q(1+b)(1+bc)+z^4a^2b^2c^2q}
%{(1-z^2ac)(1-z^2q)(1-z^2acq)(1-z^2q^2)},
\end{align*}
where $q=abcd$ as before.
If one compares these with the first few terms of $\Psi_{N}$,
one can easily guess the following theorem holds:

%----------------------------------------------------------
% Main Theorem
%----------------------------------------------------------
\begin{theorem}
\label{ordinary_strict}
For non-negative integer $N$, let $\Phi_{N}=\Phi_{N}(a,b,c,d;z)$
be as in \thetag{\ref{eq:gen_ord}} and $q=abcd$.
Then we have
\begin{align}
\Phi_{N}(a,b,c,d;z)
=\frac{\Psi_{N}(a,b,c,d;z)}{(z^2q;q)_{\lfloor N/2\rfloor}(z^2ac;q)_{\lceil N/2\rceil}},
\end{align}
where $\Psi_{N}=\Psi_{N}(a,b,c,d;z)$ is the generating function defined in \thetag{\ref{eq:gen_strict}}.
Note that $\Psi_{N}$ is explicitly given
in terms of basic hypergeometric functions
in Theorem~\ref{th:gen_strict}.
\end{theorem}

In fact,
the main purpose of this section is to prove this theorem.
Here
we give two proofs,
i.e. an algebraic proof (see Proposition~\ref{prop:rec_combi} and Proposition~\ref{prop:rec_numer})
and a bijective proof (see Remark~\ref{rem:bijection}).
Before we proceed to the proofs of this theorem
we state the corollaries immediately obtained from this theorem and the results in Section~\ref{Strict_Partitions}.
First of all,
as an immediate corollary of Theorem~\ref{ordinary_strict} and Corollary~\ref{cor:limit_strict},
we obtain the following generalization of Boulet's result
(Corollary~\ref{boulet:ordinary}).
%----------------------------------------------------------
% Corollay (limit)
%----------------------------------------------------------
\begin{corollary}
Assume $|a|,|b|,|c|,|d|<1$ and
set $q=abcd$.
Let $s^{X}_{i}$, $s^{Y}_{i}$, $X_{i}$, $Y_{i}$ ($i=0,1$) be as in Theorem~\ref{th:gen_strict}.
Then we have
\begin{align}
\sum_{\lambda}\omega(\lambda)z^{|\mu|}
&=
\frac{(-abc , -a z^2 q;q)_{\infty}}
{(ab,acz^2,z^2q;q)_{\infty}}(s^{X}_{0}X_{1} - s^{X}_{1}X_{0})
\nonumber\\
&=
\frac{( -abc , -a^2bcd z^2 q;q)_{\infty}}
{(ab,acz^2,z^2q;q)_{\infty}}
(s^{Y}_{0}Y_{1} - s^{Y}_{1}Y_{0}),
\label{eq:limit2}
\end{align}
where the sum runs over all partitions $\lambda$.
\end{corollary}
Theorem~\ref{ordinary_strict} and Theorem~\ref{hyper:strict}
also give the following corollary:
%----------------------------------------------------------
% Main Theorem
%----------------------------------------------------------
\begin{corollary}
\label{hyper:ordinary}
Put $x=\frac{(ab)^{\frac12}+(ab)^{-\frac12}}2$ and $q=abcd$.
Let $\Phi_{N}=\Phi_{N}(a,b,c,d;z)$ be as in \thetag{\ref{eq:gen_ord}}.
\begin{enumerate}
\item[(i)]
The generating function $\Phi_{2N}(a,b,c,d;1)$ is given by
\begin{align}
\Phi_{2N}(a,b,c,d;1)
&=\frac{
(ab)^{\frac{N}2}Q_{N}(x;-a^{\frac12}b^{\frac12}c,-a^{\frac12}b^{-\frac12}|q)
}
{(q;q)_{N}(ac;q)_{N}}
\nonumber\\
&=\frac{(-a;q)_{N}}{(q;q)_{N}(ac;q)_{N}}\,
\PHI{q^{-N}}{-c}{-a^{-1}q^{-N+1}}{-bq}.
\label{eq:od_al1}
\end{align}
\item[(ii)]
The generating function $\Phi_{2N}(a,b,c,d;1)$ is given by
\begin{align}
\Phi_{2N+1}(a,b,c,d;1)
&=\frac
{(1+a)(ab)^{\frac{N}2}
Q_{N}(x;-a^{\frac12}b^{\frac12}c,-a^{\frac32}b^{\frac12}cd|q)}
{(q;q)_{N}(ac;q)_{N+1}}
\nonumber\\
&=\frac{(-a;q)_{N+1}}{(q;q)_{N}(ac;q)_{N+1}}\,
\PHI{q^{-N}}{-c}{-a^{-1}q^{-N}}{-b}.
\label{eq:od_al2}
\end{align}
\end{enumerate}
\end{corollary}

Let $S_{N}(n,r,s)$ denote the number of partitions $\pi$ of $n$
where each part of $\pi$ is $\leq N$, ${\cal O}(\pi)=r$, ${\cal O}(\pi')=s$.
As before we immediately deduce the following result of Andrews (cf. \cite[Theorem~1]{And})
from Corollary~\ref{hyper:ordinary}.

%----------------------------------------------------------
% Corollary
%----------------------------------------------------------
\begin{corollary}
\label{cor:andrews}
(Andrews)
\begin{equation}
\sum_{n,r,s\geq0}S_{2N}(n,r,s)q^{n}z^{r}y^{s}
=\frac{
\sum_{j=0}^{N}
\left[{{N}\atop{j}}\right]_{q^4}
(-zyq;q^4)_{j}(-zy^{-1}q;q^4)_{N-j}(yq)^{2N-2j}
}
{(q^4;q^4)_{N}(z^2q^4;q^4)_{N}},
\end{equation}
and
\begin{equation}
\sum_{n,r,s\geq0}S_{2N+1}(n,r,s)q^{n}z^{r}y^{s}
=\frac{
\sum_{j=0}^{N}
\left[{{N}\atop{j}}\right]_{q^4}
(-zyq;q^4)_{j+1}(-zy^{-1}q;q^4)_{N-j}(yq)^{2N-2j}
}
{(q^4;q^4)_{N}(z^2q^4;q^4)_{N+1}}.
\end{equation}
\end{corollary}

Similarly, as in the strict case, we obtain immediately Boulet's corresponding result
for ordinary partitions (cf. \cite[Theorem~1]{Bou}).
%----------------------------------------------------------
% Theorem (Boulet) ordinary partitions
%----------------------------------------------------------
\begin{corollary}
\label{boulet:ordinary}
(Boulet) Let $q=abcd$, then
\begin{align}
\sum_{\lambda}\omega(\lambda)
=\frac{(-a;q)_{\infty}(-abc;q)_{\infty}}
{(q;q)_{\infty}(ab;q)_{\infty}(ac;q)_{\infty}},
\label{eq:boulet}
\end{align}
where the sum runs over all partitions.
\end{corollary}

In order to prove Theorem~\ref{ordinary_strict}
we first derive a recurrence formula for
$\Phi_{N}(a,b,c,d;z)$.
%----------------------------------------------------------
% Proposition
%----------------------------------------------------------
\begin{prop}
\label{prop:rec_combi}
Let $\Phi_{N}=\Phi_{N}(a,b,c,d;z)$ be as before and $q=abcd$.
Then the following recurrences hold for any positive integer $N$.
\begin{align}
&(1-z^2q^{N})\Phi_{2N}=(1+b)\Phi_{2N-1}-b\Phi_{2N-2},\label{eq:ord1}\\
&(1-z^2acq^{N})\Phi_{2N+1}=(1+a)\Phi_{2N}-a\Phi_{2N-1}.\label{eq:ord2}
\end{align}
\end{prop}
\begin{demo}{Proof}
It suffices to prove that
\begin{align}
&\Phi_{2N}=\Phi_{2N-1}+b(\Phi_{2N-1}-\Phi_{2N-2})+z^2q^{N}\Phi_{2N},\label{eq:ord1c}\\
&\Phi_{2N+1}=\Phi_{2N}+a(\Phi_{2N}-\Phi_{2N-1})+z^2acq^{N}\Phi_{2N+1}.\label{eq:ord2c}
\end{align}
Let ${\cal L}_{N}$ denote the set of partitions $\lambda$ such that $\lambda_1\leq N$.
The generating function of ${\cal L}_{N}$ with weight $\omega(\lambda)z^{\ell(\lambda)}$
is $\Phi_{N}=\Phi_{N}(a,b,c,d;z)$.
We divide ${\cal L}_{N}$ into three disjoint subsets:
\[
{\cal L}_{N}={\cal L}_{N-1}\uplus{\cal M}_{N}\uplus{\cal N}_{N}
\]
where ${\cal M}_{N}$ denote the set of partitions $\lambda$ such that $\lambda_1=N$ and $\lambda_2<N$,
and ${\cal N}_{N}$ denote the set of partitions $\lambda$ such that $\lambda_1=\lambda_2=N$.
When $N=2r$ is even,
it is easy to see that the generating function of ${\cal M}_{2r}$ equals $b(\Phi_{2r-1}-\Phi_{2r-2})$,
and the generating function of ${\cal N}_{2r}$ equals $z^2q^r\Phi_{2r}$.
This proves \thetag{\ref{eq:ord1c}}.
When $N=2r+1$ is odd, the same division proves \thetag{\ref{eq:ord2c}}.
$\Box$
\end{demo}

By simple computation,
one can derive the following identities from \thetag{\ref{eq:ord1}} and \thetag{\ref{eq:ord2}}.
%----------------------------------------------------------
% Proposition
%----------------------------------------------------------
\begin{prop}
\label{prop:rec_numer}
If we put
\begin{equation}\label{eq:os}
\Phi_{N}(a,b,c,d;z)
=\frac{F_{N}(a,b,c,d;z)}
{(z^2q;q)_{\lfloor{N}/2\rfloor}(z^2ac;q)_{\lceil N/2\rceil}},
\end{equation}
then,
\begin{align}
&F_{2N}=(1+b)F_{2N-1}-b(1-z^2acq^{N-1})F_{2N-2},\label{eq:strict1}\\
&F_{2N+1}=(1+a)F_{2N}-a(1-z^2q^{N})F_{2N-1}.\label{eq:strict2}
\end{align}
hold for any positive integer $N$.
\end{prop}
\begin{demo}{Proof}
Substitute \thetag{\ref{eq:os}} into \thetag{\ref{eq:ord1}} and \thetag{\ref{eq:ord2}},
and compute directly to obtain \thetag{\ref{eq:strict1}} and \thetag{\ref{eq:strict2}}.
$\Box$
\end{demo}

\begin{demo}{Proof of Theorem~\ref{ordinary_strict}}
From \thetag{\ref{eq:strict1}} and \thetag{\ref{eq:strict2}},
one easily sees that $F_{2N}(a,b,c,d;z)$ and $F_{2N+1}(a,b,c,d;z)$ satisfy exactly the same recurrence
in Theorem~\ref{th:aac}.
Further,
from the above example,
we see
\begin{align*}
&F_{0}=1,\\
&F_{1}=1+az,\\
&F_{2}=1+a(1+b)z+abcz^2,\\
&F_{3}=1+a(1+b+ab)z+abc(1+a+ad)z^2+a^3bcdz^3,\\
&F_{4}=1+a(1+b)(1+ab)z+abc(1+a+ab+ad+abd+abcd)z^2\\
&\qquad\qquad+a^3bcd(1+b)(1+bc)z^3+a^3b^3c^3dz^4.
\end{align*}
Thus the first few terms of $F_{N}(a,b,c,d;z)$ agree with those of $\Psi_{N}(a,b,c,d;z)$.
We immediately conclude that $F_{N}(a,b,c,d;z)=\Psi_{N}(a,b,c,d;z)$ for all $N$.
$\Box$
\end{demo}
%----------------------------------------------------------
% Remark
%----------------------------------------------------------
\begin{remark}
\label{rem:bijection}
Here we also give another proof of Theorem~\ref{ordinary_strict} by a bijection, which
 has already been used by Boulet~\cite{Bou} in the infinite case.
\end{remark}
\begin{demo}{Bijective proof of Theorem~\ref{ordinary_strict}}
Let ${\cal P}_N$ (resp. ${\cal D}_N$) denote the set of partitions
(resp. strict partitions) whose parts are less than or equal to
$N$ and let ${\cal E}_N$ denote the set of partitions whose parts
appear an even number of times and are less than or equal to $N$.
We shall establish a bijection $g:{\cal P}_N\longrightarrow {\cal
D}_N\times {\cal E}_N$ with $g(\lambda)=(\mu,\nu)$ defined as
follows. Suppose $\lambda$ has $k$ parts equal to $i$. If $k$ is
even then $\nu$ has $k$ parts equal to $i$, and if $k$ is odd then
$\nu$ has $k-1$ parts equal to $i$. The parts of $\lambda$ which
were not removed to form $\nu$, at most one of each cardinality,
give $\mu$. It is clear that under this bijection,
$\omega(\lambda)=\omega(\mu)\omega(\nu)$. It is easy to see that
the generating
function of ${\cal E}_N$ is equal to
$$
\prod_{j=1}^{\lfloor \frac{N}{2}\rfloor}\frac{1}{1-z^2q^j}\times
\prod_{j=0}^{\lfloor \frac{N-1}{2}\rfloor}\frac{1}{1-z^2acq^j},
$$
where $q=abcd$. As $\lfloor \frac{N-1}{2}\rfloor=\lceil \frac{N}{2}\rceil-1$, we obtain \eqref{eq:os}.
$\Box$
\end{demo}
%***% %----------------------------------------------------------
%***% % Proposition
%***% %----------------------------------------------------------
%***% \begin{prop}
%***% If we put $x_{n}=F_{2n}$ and $y_{n}=F_{2n+1}$,
%***% then we obtain the following three-term recurrences:
%***% \begin{align}
%***% &x_{n+1}
%***% =\{1+ab+z^2a(1+bc)q^{n}\}x_{n}
%***% -ab(1-z^2q^{n})(1-z^2acq^{n-1})x_{n-1},
%***% \label{eq:recx}\\
%***% &y_{n+1}
%***% =\{1+ab+z^2abc(1+ad)q^{n}\}y_{n}
%***% -ab(1-z^2q^{n})(1-z^2a^2bc^2dq^{n-1})y_{n-1},
%***% \label{eq:recy}
%***% \end{align}
%***% with
%***% $x_0=1$, $x_{1}=1+za(1+b)+z^2abc$, $y_0=1+za$
%***% and
%***% \[
%***% y_1=1+za(1+b+ab)+z^2a(bc+abc+q)+z^3a^2q.
%***% \]
%***% \end{prop}
%***% \begin{demo}{Proof of Theorem~\ref{hyper:ordinary}}
%***% Note that \thetag{\ref{eq:recx}} and \thetag{\ref{eq:recy}}
%***% are exactly the same as \thetag{\ref{eq:recX}} and \thetag{\ref{eq:recY}},
%***% respectively.
%***% We don't have any explanation on this coincidence at this point.
%***% But, especially,
%***% if $z=1$,
%***% then the initial condition are also the same.
%***% This leads to the conclusion.
%***% $\Box$
%***% \end{demo}
%***%
\bigskip

At the end of this section we state another enumeration of the ordinary partitions,
which is not directly related to Andrews' result,
but obtained as an application of the minor summation formula of Pfaffians.
%Our goal is to prove the following theorem.
Let
\[
\Phi_{N,M}
=\Phi_{N,M}(a,b,c,d)
=\sum_{{\lambda}\atop{\lambda_1\leq N,\ \ell(\lambda)\leq M}}\omega(\lambda),
\]
where the sum runs over all partitions $\lambda$
such that $\lambda$ has at most $M$ parts and each part of $\lambda$ is less than or equal to $N$.

Again we use Lemma~\ref{lem:pf_sum} and Theorem~\ref{thm:coefficients}
to obtain the following theorem.
\begin{theorem}
Let $N$ be a positive integer and set $q=abcd$.
Then we have
\begin{align}
&\sum_{t=0}^{\lfloor N/2\rfloor}
\Phi_{N-2t,2t}(a,b,c,d)\,
z^t q^{\binom{t}2}
%\sum_{{\lambda} \atop {{\ell(\lambda)\leq 2t} \atop {\lambda_1\leq N-2t}}}
%\omega(\lambda)
=\Pf\begin{bmatrix}
S&J\\
-J&C
\end{bmatrix},
\label{eq:pf_ord}
\end{align}
where $S=(1)_{0\leq i<j\leq N-1}$ and
$C=(
a^{\lceil (j-1)/2\rceil}
b^{\lfloor (j-1)/2\rfloor}
c^{\lceil i/2\rceil}
d^{\lfloor i/2\rfloor}z
)_{0\leq i<j\leq N-1}$.
\end{theorem}
\begin{demo}{Proof}
As before,
we take $A=(1)_{0\leq i<j\leq N-1}$ and
$$
B=(a^{\lceil (j-1)/2\rceil}b^{\lfloor (j-1)/2\rfloor}
c^{\lceil i/2\rceil}d^{\lfloor i/2\rfloor})_{0\leq i<j\leq N-1},
$$
in Lemma~\ref{lem:pf_sum},
then \thetag{\ref{eq:pf_ord}} follows from Lemma~\ref{thm:coefficients}.
$\Box$
\end{demo}
For example,
if $N=4$,
then  the right-hand side of \thetag{\ref{eq:pf_ord}} becomes
\[
\Pf\left[ \begin {array}{cccc|cccc}
0&1&1&1&0&0&0&1\\
-1&0&1&1&0&0&1&0\\
-1&-1&0&1&0&1&0&0\\
-1&-1&-1&0&1&0&0&0\\\hline
0&0&0&-1&0&z&az&abz\\
0&0&-1&0&-z&0&acz&abcz\\
0&-1&0&0&-az&-acz&0&abcdz\\
-1&0&0&0&-abz&-abcz&-abcdz&0
\end {array} \right].
\]
Let $\tilde\Phi_N=\tilde\Phi_N(a,b,c,d;z)=\Pf\begin{bmatrix}
S&J\\
-J&C
\end{bmatrix}$
denote the right-hand side of \thetag{\ref{eq:pf_ord}}.
For example,
we have
$\tilde\Phi_1=1$,
$\tilde\Phi_2=1+z$,
$\tilde\Phi_3=1+(1+a+ac)z$
and
$\tilde\Phi_4=1+(1+a+ab+ac+abc+abcd)z+abcdz^2$.
Note that the partitions $\lambda$
such that $\ell(\lambda)\leq2$ and $\lambda_1\leq2$ are the following six:
\[
\emptyset
\qquad
\smyoung{
a\\
}
\qquad
\smyoung{
a&b\\
}
\qquad
\smyoung{
a\\
c\\
}
\qquad
\smyoung{
a&b\\
c\\
}
\qquad
\smyoung{
a&b\\
c&d\\
}\ .
\]
The sum of their weights is equal to $[z]\tilde\Phi_4=1+a+ab+ac+abc+abcd$.

The same argument as in the proof of Proposition~\ref{prop:rec_dist}
can be used to prove the following proposition.
%----------------------------------------------------------
% Proposition
%----------------------------------------------------------
\begin{prop}
\label{prop:rec_ord}
Let $\tilde\Phi_N=\tilde\Phi_N(a,b,c,d;z)$ be as above.
Then we have
\begin{align}
&\tilde\Phi_{2N}=(1+b)\tilde\Phi_{2N-1}+(a^{N-1}b^{N-1}c^{N-1}d^{N-1}z-b)\tilde\Phi_{2N-2},
\label{eq:o1}\\
&\tilde\Phi_{2N+1}=(1+a)\tilde\Phi_{2N}+(a^{N}b^{N-1}c^{N}d^{N-1}z-a)\tilde\Phi_{2N-1},
\label{eq:o2}
\end{align}
for any positive integer $N$.
\end{prop}
\begin{demo}{Proof}
Perform the same elementary transformations of rows and columns on $\begin{bmatrix}
S&J\\
-J&C
\end{bmatrix}$
as we did in the proof of Proposition~\ref{prop:rec_dist},
and expand it along the last row/column.
The details are left to the reader.
$\Box$
\end{demo}
\begin{remark}
The recurrence equations \thetag{\ref{eq:o1}} and \thetag{\ref{eq:o2}} also can be proved combinatorially.
\end{remark}
\begin{demo}{Proof of Proposition~\ref{prop:rec_ord}}
Consider the generating function of partitions:
\begin{equation}\label{eq:p1}
\sum_{{\lambda} \atop {{\ell(\lambda)\leq 2t} \atop {\lambda_1\leq 2j+1-2t}}}w(\lambda)=
\sum_{{\lambda} \atop {{\ell(\lambda)\leq 2t} \atop {\lambda_1\leq 2j-2t}}}w(\lambda)+
\sum_{{\lambda} \atop {{\ell(\lambda)\leq 2t} \atop {\lambda_1=2j+1-2t}}}w(\lambda).
\end{equation}
Splitting the  partitions $\lambda$ in the  second sum of the right side
into two subsets:
$\lambda_2<\lambda_1$, and $\lambda_2=\lambda_1$. Now
\begin{equation}\label{eq:p2}
\sum_{{\lambda: \lambda_1>\lambda_2} \atop {{\ell(\lambda)\leq 2t} \atop {\lambda_1=2j+1-2t}}}w(\lambda)=
a\left(\sum_{{\lambda} \atop {{\ell(\lambda)\leq 2t} \atop {\lambda_1\leq 2j-2t}}}w(\lambda)-
\sum_{{\lambda} \atop {{\ell(\lambda)\leq 2t} \atop {\lambda_1\leq 2j-1-2t}}}w(\lambda)\right),
\end{equation}
and
\begin{equation}\label{eq:p3}
\sum_{{\lambda: \lambda_1=\lambda_2} \atop {{\ell(\lambda)\leq 2t} \atop {\lambda_1=2j+1-2t}}}w(\lambda)=
acq^{j-t}\sum_{{\lambda} \atop {{\ell(\lambda)\leq 2t-2} \atop {\lambda_1\leq 2j+1-2t}}}w(\lambda).
\end{equation}
 Plugging \eqref{eq:p2} and \eqref{eq:p3} into \eqref{eq:p1} and then multiplying
 by $z^tq^{t\choose 2}$ and summing over $t$ we get \eqref{eq:o2}.
 Similarly we can prove \eqref{eq:o1}.
$\Box$
\end{demo}
%----------------------------------------------------------
% Proposition
%----------------------------------------------------------
\begin{prop}
 Set $U_N=\tilde\Phi_{2N}$ and $V_N=\tilde\Phi_{2N+1}$, then, for $N\geq 1$,
\begin{align}
U_{N+1}&= \left\{1+ab+ac(1+bd)q^{N-1}z\right\}U_{N}-
a(b-zq^{N-1})(1-czq^{N-1})U_{N-1},\\
V_{N+1}&= \left\{1+ab+(1+ac)zq^{N}\right\}V_N-a(b-zq^{N})(1-czq^{N-1})
V_{N-1},
\end{align}
where
$U_0=1$, $V_0=1$, $U_1=1+z$, $V_1=1+(1+a+ac)z$.
\end{prop}

Thus $U_N$ and $V_N$ are also expressed by the solutions of the associated Al-Salam-Chihara polynomials.

\bigskip

\section{A weighted sum of Schur's $P$-functions}

We use the notation $X=X_n=(x_1,\dots,x_n)$ for the finite set of variables $x_1$, $\dots$, $x_n$.
%In \cite{I},
%one of the authors used a Pfaffian expression of
%$\displaystyle
%\sum_{\lambda}\omega(\lambda) s_{\lambda}(X)$
%to prove Stanley's open problem,
%where the sum runs over all partition $\lambda$
%and $s_{\lambda}(X)$ stands for the Schur function with respect to a partition $\lambda$.
%***% Let us call $P_{\mu}(X)=P_{\mu}(X,-1)$ \defterm{Schur's $P$-function},
%***% and $Q_{\mu}(X)=Q_{\mu}(X,-1)$ \defterm{Schur's $Q$-function}
%***% where $P_{\mu}(X,t)$ and $Q_{\mu}(X,t)$ are Hall-Littlewood polynomials.
The aim of this section is to give some Pfaffian and determinantal formulas for the weighted sum
$\sum \omega(\mu)z^{\ell(\mu)}P_{\mu}(x)$ where $P_{\mu}(x)$ is Schur's $P$-function.

Let $A_n$ denote the skew-symmetric matrix
\begin{align*}
\left(\frac{x_i-x_j}{x_i+x_j}\right)_{1\leq i,j\leq n}
\end{align*}
and for each strict partition $\mu=(\mu_1,\dots,\mu_l)$ of length $l\leq n$,
let $\Gamma_{\mu}$ denote the $n\times l$ matrix $\left(x_j^{\mu_i}\right)$.
Let
\begin{equation*}
A_{\mu}(x_1,\dots,x_n)=\begin{pmatrix}
A_n&\Gamma_{\mu}J_{l}\\
-J_l{}^t\!\Gamma_{\mu}&O_l
\end{pmatrix}
\end{equation*}
which is a skew-symmetric matrix of $(n+l)$ rows and columns.
Define $\Pf_\mu(x_1,\dots,x_n)$ to be $\Pf A_{\mu}(x_1,\dots,x_n)$
if $n+l$ is even,
and to be $\Pf A_{\mu}(x_1,\dots,x_n,0)$
if $n+l$ is odd.
By \cite[Ex.13, p.267]{Ma},
Schur's $P$-function $P_\mu(x_1,\dots,x_n)$ is defined to be
\begin{equation*}
\frac{\Pf_\mu(x_1,\dots,x_n)}{\Pf_\emptyset(x_1,\dots,x_n)},
\end{equation*}
where it is well-known that
$
\Pf_\emptyset(x_1,\dots,x_n)
=\prod_{1\leq i<j\leq n}
\frac{x_i-x_j}{x_i+x_j}.
$
Meanwhile,
by \cite[(8.7), p.253]{Ma},
Schur's $Q$-function $Q_\mu(x_1,\dots,x_n)$ is defined to be
$
2^{\ell(\lambda)}P_\mu(x_1,\dots,x_n)
$.
%These are not ordinary way to define the $P$-functions and $Q$-functions,
%but here we adopt these definitions
%since it makes us easier to treat Pfaffians.

In this section,
we consider a weighted sum of Schur's $P$-functions and $Q$-functions,
i.e.,
\begin{align*}
&\xi_{N}(a,b,c,d;X_n)
=\sum_{{\mu}
\atop{\mu_1\leq N}}
\omega(\mu)P_{\mu}(x_1,\dots,x_n),
\\
&\eta_{N}(a,b,c,d;X_n)
=\sum_{{\mu}
\atop{\mu_1\leq N}}
\omega(\mu)Q_{\mu}(x_1,\dots,x_n),
\end{align*}
where the sums run over all strict partitions $\mu$
such that each part of $\mu$ is less than or equal to $N$.
More generally,
we can unify these problems to finding the following sum:
\begin{equation}
\zeta_{N}(a,b,c,d;z;X_n)
=\sum_{{\mu}
\atop{\mu_1\leq N}}
\omega(\mu)
z^{\ell(\mu)}
P_{\mu}(x_1,\dots,x_n),
\end{equation}
where the sum runs over all strict partitions $\mu$
such that each part of $\mu$ is less than or equal to $N$.
One of the main results of this section is that $\zeta_{N}(a,b,c,d;z;X_n)$
can be expressed by a Pfaffian
(see Corollary~\ref{cor:sumQ_finite}).
Further,
let us put
\begin{align}
\zeta(a,b,c,d;z;X_n)=\lim_{N\rightarrow\infty}\zeta_{N}(a,b,c,d;z;X_n)
=\sum_{\mu}
\omega(\mu)
z^{\ell(\mu)}
P_{\mu}(X_n),
\end{align}
where the sum runs over all strict partitions $\mu$.
We also write
\begin{align*}
\xi(a,b,c,d;X_n)=\zeta(a,b,c,d;1;X_n)
=\sum_{\mu}
\omega(\mu)
P_{\mu}(X_n),
\end{align*}
where the sum runs over all strict partitions $\mu$.
Then we have the following theorem:
%----------------------------------------------------------
% Corollary
%----------------------------------------------------------
\begin{theorem}
\label{cor:sumQ_infinite}
Let $n$  be a positive integer.
Then
\begin{align}
\zeta(a,b,c,d;z;X_{n})
=\begin{cases}
\Pf\left(\gamma_{ij}\right)_{1\leq i<j\leq n}/\Pf_\emptyset(X_n)
&\text{ if $n$ is even,}\\
\Pf\left(\gamma_{ij}\right)_{0\leq i<j\leq n}/\Pf_\emptyset(X_n)
&\text{ if $n$ is odd,}
\end{cases}
\label{eq:zeta_PF}
\end{align}
where
\begin{align}
\gamma_{ij}
&=\frac{x_{i}-x_{j}}{x_{i}+x_{j}}
+u_{ij}z+v_{ij}z^2
\end{align}
with
\begin{align}
u_{ij}
&=\frac{a\det\begin{pmatrix}x_{i}+bx_{i}^{2}&1-abx_{i}^{2}\\ x_{j}+bx_{j}^{2}&1-abx_{j}^{2}\end{pmatrix}}
{(1-abx_i^2)(1-abx_j^2)},
\\
v_{ij}
&=\frac{abcx_{i}x_{j}\det\begin{pmatrix}x_{i}+ax_{i}^{2}&1-a(b+d)x_{i}^{2}-abdx_{i}^{3}\\ x_{j}+ax_{j}^{2}&1-a(b+d)x_{j}^{2}-abdx_{j}^{3}\end{pmatrix}}
{(1-abx_i^2)(1-abx_j^2)(1-abcdx_i^2x_j^2)},
\end{align}
if $1\leq i,j\leq n$,
and
\begin{align}
\gamma_{0j}
=1+\frac{ax_{j}(1+bx_{j})}{1-abx_{j}^2}z
\end{align}
if $1\leq j\leq n$.

Especially, when $z=1$, we have
\begin{align}
\xi(a,b,c,d;X_{n})
=\begin{cases}
\Pf\left({\widetilde\gamma}_{ij}\right)_{1\leq i<j\leq n}/\Pf_\emptyset(X_n)
&\text{ if $n$ is even,}\\
\Pf\left({\widetilde\gamma}_{ij}\right)_{0\leq i<j\leq n}/\Pf_\emptyset(X_n)
&\text{ if $n$ is odd,}
\end{cases}
\end{align}
where
\begin{align}
{\widetilde\gamma}_{ij}
=\begin{cases}
\frac{1+ax_{j}}{1-abx_{j}^2}
&\text{ if $i=0$,}\\
\frac{x_{i}-x_{j}}{x_{i}+x_{j}}
+\widetilde{v}_{ij}
&\text{ if $1\leq i<j\leq n$,}
\end{cases}
with
\end{align}
\begin{equation}
\widetilde{v}_{ij}=\frac{a\det\begin{pmatrix}x_{i}+bx_{i}^{2}&1-b(a+c)x_{i}^{2}-abcx_{i}^{3}\\ x_{j}+bx_{j}^{2}&1-b(a+c)x_{j}^{2}-abcx_{j}^{3}\end{pmatrix}}
{(1-abx_i^2)(1-abx_j^2)(1-abcdx_i^2x_j^2)}.
\end{equation}
\end{theorem}

We can generalize
this result in the following theorem (Theorem~\ref{th:xi})
using the generalized Vandermonde determinant used in \cite{IOTZ}.
Let $n$ be an non-negative integer,
and let $X=(x_{1},\dots,x_{2n})$, $Y=(y_{1},\dots,y_{2n})$,
$A=(a_{1},\dots,a_{2n})$ and $B=(b_{1},\dots,b_{2n})$ be $2n$-tuples of variables.
Let $V^{n}(X,Y,A)$ denote the $2n\times n$ matrix whose $(i,j)$th entry is
$a_ix_{i}^{n-j}y_{i}^{j-1}$ for $1\leq i\leq 2n$, $1\leq j\leq n$,
and let $U^n(X,Y;A,B)$ denote the $2n\times 2n$ matrix
$
\begin{pmatrix}V^{n}(X,Y,A)&V^{n}(X,Y,B)\end{pmatrix}.
$
For instance if $n=2$ then $U^2(X,Y;A,B)$ is
\[
\begin{pmatrix}
a_{1}x_{1}&a_{1}y_{1}&b_{1}x_{1}&b_{1}y_{1}\\
a_{2}x_{2}&a_{2}y_{2}&b_{2}x_{2}&b_{2}y_{2}\\
a_{3}x_{3}&a_{3}y_{3}&b_{3}x_{3}&b_{3}y_{3}\\
a_{4}x_{4}&a_{4}y_{4}&b_{4}x_{4}&b_{4}y_{4}
\end{pmatrix}.
\]
Hereafter we use the following notation for
 $n$-tuples $X = (x_1, \cdots, x_n)$ and $Y = (y_1, \cdots, y_n)$ of variables:
$$
X + Y
 = (x_1 + y_1, \ldots, x_n + y_n),
\quad
X\cdot Y
 = (x_1 y_1, \ldots, x_n y_n),
$$
and, for integers $k$ and $l$,
$$
X^k
 = (x_1^k, \ldots, x_n^k),
\quad
X^k Y^l
 = (x_1^k y_1^l, \ldots, x_n^k y_n^l).
$$
Let $\pmb1$ denote the $n$-tuple $(1,\dots,1)$.
For any subset $I=\{i_1,\dots,i_r\}\in\binom{[n]}{r}$,
let $X_{I}$ denote the $r$-tuple $(x_{i_1},\dots,x_{i_r})$.
%----------------------------------------------------------
% Main theorem
%----------------------------------------------------------
\begin{theorem}
\label{th:xi}
%Let $n$  be a positive integer.
Let $q=abcd$.
If $n$ is an even integer,
then we have
\begin{align}
\xi(a,b,c,d;X_{n})
&=\sum_{r=0}^{n/2}
\sum_{I\in\binom{[n]}{2r}}
\frac{(-1)^{|I|-\binom{r+1}{2}}a^rq^{\binom{r}{2}}}{\prod_{i\in I}(1-abx_i^2)}
\prod_{{i,j\in I}\atop{i<j}}\frac{x_i+x_j}{(x_i-x_j)(1-qx_i^2x_j^2)}
\nonumber\\
&\times
\det U^{r}(X_I^2,\pmb1+qX_I^4,X_I+bX_I^2,\pmb1-b(a+c)X_{I}^2-abcX_{I}^3).
\label{eq:det_xi_even}
\end{align}
If $n$ is an odd integer,
then we have
\begin{align}
&\xi(a,b,c,d;X_{n})
=\sum_{m=1}^{n}\frac{1+ax_{m}}{1-abx_{m}^2}\sum_{r=0}^{(n-1)/2}
\sum_{I\in\binom{[n]\setminus\{m\}}{2r}}
\frac{(-1)^{|I|-\binom{r+1}{2}}a^rq^{\binom{r}{2}}}{\prod_{i\in I}(1-abx_i^2)}
\prod_{i\in I}\frac{x_{m}+x_{i}}{x_{m}-x_{i}}
\nonumber\\
&\times
\prod_{{i,j\in I}\atop{i<j}}\frac{x_i+x_j}{(x_i-x_j)(1-qx_i^2x_j^2)}
\cdot\det U^{r}(X_I^2,\pmb1+qX_I^4,X_I+bX_I^2,\pmb1-b(a+c)X_{I}^2-abcX_{I}^3).
\label{eq:det_xi_odd}
\end{align}
\end{theorem}
%----------------------------------------------------------
% Main theorem
%----------------------------------------------------------
\begin{theorem}
\label{th:zeta}
%Let $n$  be a positive integer.
Let $q=abcd$.
If $n$ is an even integer,
then $\zeta(a,b,c,d;z;X_{n})$ is equal to
\begin{align}
&\sum_{r=0}^{n/2}z^{2r}
\sum_{I\in\binom{[n]}{2r}}
\frac{(-1)^{|I|-\binom{r+1}{2}}(abc)^rq^{\binom{r}{2}}\prod_{i\in I}x_{i}}
{\prod_{i\in I}(1-abx_i^2)}
\prod_{{i,j\in I}\atop{i<j}}\frac{x_i+x_j}{(x_i-x_j)(1-qx_i^2x_j^2)}
\nonumber\\
&\times
\det U^{r}(X_I^2,\pmb1+qX_I^4,X_I+aX_I^2,\pmb1-a(b+d)X_{I}^2-abdX_{I}^3)
\nonumber\\
&+\sum_{r=0}^{n/2}z^{2r-1}
\sum_{I\in\binom{[n]}{2r}}
\sum_{{k<l}\atop{k,l\in I}}
\frac{(-1)^{|I|-\binom{r}{2}-1}a^rb^{r-1}c^{r-1}q^{\binom{r-1}{2}}\{1+b(x_k+x_l)+abx_kx_l\}\prod_{i\in I'}x_{i}}
{\prod_{i\in I}(1-abx_i^2)}
\nonumber\\
&\times\frac{\prod_{{i,j\in I}\atop{i<j}}(x_i+x_j)\cdot
\det U^{r-1}(X_{I'}^2,\pmb1+qX_{I'}^4,X_{I'}+aX_{I'}^2,\pmb1-a(b+d)X_{I'}^2-abdX_{I'}^3)}
{\prod_{{i,j\in I'}\atop{i<j}}(x_i-x_j)(1-qx_i^2x_j^2)},
\label{eq:det_zeta_even}
\end{align}
where $I'=I\setminus\{k,l\}$.
\end{theorem}
Note that we can obtain a similar formula when $n$ is odd
by expanding the Pfaffian in \thetag{\ref{eq:zeta_PF}} along the first row/column.

To obtain the sum of this type we need a generalization of Lemma~\ref{lem:pf_sum},
in which the row/column indices always contain say the set $\{1,2,...,n\}$, for some fixed $n$.
%----------------------------------------------------------
% Lemma (Minor Sum Formula)
%----------------------------------------------------------
\begin{lemma}
\label{lem:pf_sum2}
Let $n$ and $N$ be nonnegative integers.
Let $A=(a_{ij})$ and $B=(b_{ij})$
be skew symmetric matrices of size $(n+N)$.
We divide the set of row/column indices
into two subsets,
i.e. the first $n$ indices $I_0=[n]$ and the last $N$ indices $I_1=[n+1,n+N]$.
Then
\begin{align}
&\sum_{{t\geq0}
\atop{n+t\text{ even}}}
z^{(n+t)/2}
\sum_{I\in\binom{I_1}{t}}
\gamma^{|I_0\uplus I|}
\Pf\left(\Delta^{I_0\uplus I}_{I_0\uplus I}(A)\right)
\Pf\left(\Delta^{I_0\uplus I}_{I_0\uplus I}(B)\right)
\nonumber\\
&\qquad\qquad=\Pf\begin{pmatrix}
J_{n+N}\,{}^t\kern-2pt AJ_{n+N}&K_{n,N}\\
-{}^t\!K_{n,N}&C
\end{pmatrix},
\label{eq:pf_sum2}
\end{align}
where $C=(C_{ij})_{1\leq i,j,\leq n+N}$ is given by
$C_{ij}=\gamma^{i+j}b_{ij}z$
and $K_{n,N}=J_{n+N}{\widetilde E}_{n,N}$ with
\begin{equation*}
{\widetilde E}_{n,N}=\begin{pmatrix}
O_{n}&O_{n,N}\\
O_{N,n}&E_{N}
\end{pmatrix}.
\end{equation*}
\end{lemma}
%----------------------------------------------------------
% Proof of Lemma
%----------------------------------------------------------
\begin{demo}{Proof}
Let $V=\{(n+N)^{*},\dots,(n+1)^{*},n^{*},\dots,1^{*},1,\dots,n,n+1,\dots,n+N\}$
be vertices arranged in this order on the $x$-axis.
Put $V_0=\{n^{*},\dots,1^{*}\}$ and $V_1^{*}=\{(n+N)^{*},\dots,(n+1)^{*}\}$, $V_0=\{1,\dots,n\}$ and $V_1=\{n+1,\dots,N\}$.
From \thetag{\ref{eq:Pf_matching}},
the Pfaffian on the right-hand side of \thetag{\ref{eq:pf_sum2}} is equal to
\begin{equation*}
\sum_{\sigma}\sgn\sigma\,
\prod_{(i,j)\in\sigma}a_{ij}
\prod_{(i,j)\in\sigma}C_{ij}
\end{equation*}
summed over all perfect matching $\sigma$
on $V$,
in which
there is a set $I=\{i_1,\dots,i_t\}\in V_{1}$ of $t$ vertices
such that $\sigma|_{V_0\uplus V_1}$ is a perfect matching on $V_0\uplus I$
and $\sigma|_{V_0^*\uplus V_1^*}$ is a perfect matching on $V_0^*\uplus I^*$,
where $I^*=\{i_1^*,\dots,i_t^*\}$,
and each $j\in V_{1}\setminus I$ is adjoint to $j^*\in V_{1}^*\setminus I^*$ in $\sigma$.
Thus the summand vanishes unless $n+t$ is even,
and this sum is equal to
\begin{equation*}
\sum_{t}z^{(t+n)/2}\sum_{I\in\binom{V_1}{t}}\gamma^{n+|I|}
\sum_{(\sigma_1,\sigma_2)}\sgn\sigma_1\sgn\sigma_2\,
\prod_{(i,j)\in\sigma_1}a_{ij}
\prod_{(i,j)\in\sigma_2}b_{ij},
\end{equation*}
where the third sum runs over all pairs $(\sigma_1,\sigma_2)$ where $\sigma_1$ is a perfect matching on $V_0\uplus I$
and $\sigma_2$ is a perfect matching on $V_0^*\uplus I^*$.
This is equal to the left-hand side of \thetag{\ref{eq:pf_sum2}}.
$\Box$
\end{demo}

For a nonnegative integer $N$,
let $\mu^N=(N,\dots,1,0)$,
and let $\Gamma_{\mu^N}$ denote the $n\times (N+1)$ matrix $\left(x_{i}^{N-j}\right)_{1\leq i\leq n,0\leq j\leq N}$.
Let
\[
{\cal A}_{n,N}
=\begin{pmatrix}
A_n&\Gamma_{\mu^N}J_{N+1}\\
-J_{N+1}{}^t\!\Gamma_{\mu^N}&O_{N+1}
\end{pmatrix}
\]
which is a skew-symmetric matrix of size $n+N+1$.
For example,
if $n=4$ and $N=3$,
then
\[
{\cal A}_{4,3}=\begin{pmatrix}
0&{\frac {x_{{1}}-x_{{2}}}{x_{{1}}+x_{{2}}}}&{\frac {x_{{1}}-x_{{3}}}{x_{{1}}+x_{{3}}}}&{\frac {x_{{1}}-x_{{4}}}{x_{{1}}+x_{{4}}}}&1&x_{{1}}&{x_{{1}}}^{2}&{x_{{1}}}^{3}\\
{\frac {x_{{2}}-x_{{1}}}{x_{{1}}+x_{{2}}}}&0&{\frac {x_{{2}}-x_{{3}}}{x_{{2}}+x_{{3}}}}&{\frac {x_{{2}}-x_{{4}}}{x_{{2}}+x_{{4}}}}&1&x_{{2}}&{x_{{2}}}^{2}&{x_{{2}}}^{3}\\
{\frac {x_{{3}}-x_{{1}}}{x_{{1}}+x_{{3}}}}&{\frac{x_{{3}}-x_{{2}}}{x_{{2}}+x_{{3}}}}&0&{\frac {x_{{3}}-x_{{4}}}{x_{{3}}+x_{{4}}}}&1&x_{{3}}&{x_{{3}}}^{2}&{x_{{3}}}^{3}\\
{\frac {x_{{4}}-x_{{1}}}{x_{{1}}+x_{{4}}}}&{\frac {x_{{4}}-x_{{2}}}{x_{{2}}+x_{{4}}}}&{\frac {x_{{4}}-x_{{3}}}{x_{{3}}+x_{{4}}}}&0&1&x_{{4}}&{x_{{4}}}^{2}&{x_{{4}}}^{3}\\
-1&-1&-1&-1&0&0&0&0\\
-x_{{1}}&-x_{{2}}&-x_{{3}}&-x_{{4}}&0&0&0&0\\
-{x_{{1}}}^{2}&-{x_{{2}}}^{2}&-{x_{{3}}}^{2}&-{x_{{4}}}^{2}&0&0&0&0\\
-{x_{{1}}}^{3}&-{x_{{2}}}^{3}&-{x_{{3}}}^{3}&-{x_{{4}}}^{3}&0&0&0&0
\end{pmatrix}.
\]
Let $\beta_{ij}$ be as in \thetag{\ref{eq:strict}}.
Let $B_{N}$ denote the $(N+1)\times(N+1)$ matrix $(\beta_{ij})_{0\leq i,j\leq N}$
and let $B_{N}'$ denote the $(N+2)\times(N+2)$ matrix $(\beta_{ij})_{-1\leq i,j\leq N}$.
%----------------------------------------------------------
% Theorem
%----------------------------------------------------------
\begin{theorem}
Let $n$ and $N$ be integers such that $n\geq N\geq0$.
Then
\begin{align}
\label{eq:SchurQ_Pf}
\zeta_{N}(a,b,c,d;z;X_{n})
=\Pf\left({\cal C}_{n,N}\right)/\Pf_\emptyset(X_n),
\end{align}
where
\begin{align}
\label{eq:SchurQ_Pf_even}
{\cal C}_{n,N}
=\begin{pmatrix}
O_{N+1}&{}^t\!\Gamma_{\mu^N}J_{n}&J_{N+1}\\
-J_{n}\Gamma_{\mu^N}&J_{n}{}^t\!A_nJ_{n}&O_{n,N+1}\\
-J_{N+1}&O_{N+1,n}&B_{N}
\end{pmatrix},
\end{align}
if $n$ is even,
and
\begin{align}
\label{eq:SchurQ_Pf_odd}
{\cal C}_{n,N}
=\begin{pmatrix}
O_{N+1}&{}^t\!\Gamma_{\mu^N}J_{n}&J_{N+1}'\\
-J_{n}\Gamma_{\mu^N}&J_{n}{}^t\!A_nJ_{n}&O_{n,N+2}\\
-{}^t\!J_{N+1}'&O_{N+2,n}&B_{N}'
\end{pmatrix}
\end{align}
where $J_{N+1}'=\begin{pmatrix}O_{N+1,1}&J_{N+1}\end{pmatrix}$ if $n$ is odd.
\end{theorem}
%----------------------------------------------------------
% Proof of Theorem
%----------------------------------------------------------
\begin{demo}{Proof}
Let ${\cal B}_{n,N}$ be the skew-symmetric matrix of size $(n+N+1)$ defined by
\[
{\cal B}_{n,N}
=\begin{pmatrix}
S_{n}&O_{n,N+1}\\
O_{N+1.n}&B_{N}
\end{pmatrix}
\]
if $n$ is even, and
\[
{\cal B}_{n,N}
=\begin{pmatrix}
S_{n-1}&O_{n,N+2}\\
O_{N+2.n}&B_{N}'
\end{pmatrix}
\]
if $n$ is odd.
Fix a strict partition $\mu=(\mu_1,\dots,\mu_{l})$
such that $\mu_1>\dots>\mu_l\geq0$,
and let $K_n(\mu)=\{n+\mu_{l},\dots,n+\mu_{1}\}$.
From the definition of ${\cal B}_{n,N}$ and Theorem~\ref{thm:coef_strict},
we have
\begin{align*}
\Pf\left(\Delta_{[n]\uplus K_n(\mu)}^{[n]\uplus K_n(\mu)}
\left({\cal B}_{n,N}\right)\right)
=\omega(\mu)\, z^{\ell(\mu)}
\end{align*}
if $n+l$ is even.
Thus Lemma~\ref{lem:pf_sum2} immediately implies that
$\Pf_\emptyset(X_n)\zeta_{N}(a,b,c,d;z;X_{n})$ is equal to
\begin{equation}
\Pf\begin{pmatrix}
J_{n+N+1}{}^t\!{\cal A}_{n,N}J_{n+N+1}&K_{n,N+1}\\
-{}^{t}\!K_{n,N+1}&{\cal B}_{n,N}
\end{pmatrix}.
\end{equation}
By simple elementary transformations on rows and columns,
we obtain the desired results \thetag{\ref{eq:SchurQ_Pf_even}} and \thetag{\ref{eq:SchurQ_Pf_odd}}.
$\Box$
\end{demo}
For instance,
if $n=4$ and $N=2$,
then ${\cal D}_{4,2}$ looks as follows:
\[
\begin{pmatrix}
 0& 0& 0&x_4^2&x_3^2&x_2^2&x_1^2& 0& 0& 1\\
 0& 0& 0&x_4&x_3&x_2&x_1& 0& 1& 0\\
 0& 0& 0& 1& 1& 1& 1& 1& 0& 0\\
-x_4^2&-x_4&-1&0&\frac{x_{3}-x_{4}}{x_{3}+x_{4}}&\frac{x_{2}-x_{4}}{x_{2}+x_{4}}&\frac{x_{1}-x_{4}}{x_{1}+x_{4}}& 0& 0& 0\\
-x_3^2&-x_3&-1&\frac{x_{4}-x_{3}}{x_{4}+x_{3}}& 0& \frac {x_{{2}}-x_{{3}}}{x_{{2}}+x_{{3}}}& \frac {x_{{1}}-x_{{3}}}{x_{{1}}+x_{{3}}}& 0& 0& 0\\
-x_2^2&-x_2&-1&\frac{x_{4}-x_{2}}{x_{4}+x_{2}}&\frac{x_{3}-x_{2}}{x_{3}+x_{2}}& 0& \frac {x_{{1}}-x_{{2}}}{x_{{1}}+x_{{2}}}& 0& 0& 0\\
-x_1^2&-x_1&-1&\frac{x_{4}-x_{1}}{x_{4}+x_{1}}&\frac{x_{3}-x_{1}}{x_{3}+x_{1}}&\frac{x_{2}-x_{1}}{x_{2}+x_{1}}& 0& 0& 0& 0\\
 0& 0&-1& 0& 0& 0& 0&   0&   az&  abz\\
 0&-1& 0& 0& 0& 0& 0&  -az&   0& abcz^2\\
-1& 0& 0& 0& 0& 0& 0& -abz&-abcz^2&   0\\
\end{pmatrix}
\]
%----------------------------------------------------------
% Corollary
%----------------------------------------------------------
\begin{corollary}
\label{cor:sumQ_finite}
Let $n$ and $N$ be integers such that $n\geq N\geq0$.
Then
\begin{align}
\zeta_{N}(a,b,c,d;z;X_{n})
=\Pf\left({\cal D}_{n,N}\right)/\Pf_\emptyset(X_n),
\end{align}
where
\begin{align}
{\cal D}_{n,N}
=\left(
\frac{x_{i}-x_{j}}{x_{i}+x_{j}}
+\sum_{0\leq k,l\leq N}\beta_{kl}x_{i}^{l}x_{j}^{k}
%\begin{vmatrix}x_{i}^{k}&x_{i}^{l}\\ x_{j}^{k}&x_{j}^{l}\end{vmatrix}
\right)_{1\leq i,j\leq n},
\end{align}
if $n$ is even,
and
\begin{align}
{\cal D}_{n,N}
=\left(
\begin{array}{c|c}
0&{\displaystyle\sum_{k=0}^{N}\beta_{-1,k}x_j^k}\\
\hline
{\displaystyle\sum_{k=0}^{N}\beta_{k,-1}x_i^k}&
{\displaystyle
\frac{x_{i}-x_{j}}{x_{i}+x_{j}}
+\sum_{0\leq k,l\leq N}\beta_{kl}x_{i}^{l}x_{j}^{k}}
\end{array}
\right)_{0\leq i,j\leq n},
\end{align}
if $n$ is odd.
\end{corollary}
\begin{demo}{Proof}
When $n$ is even,
annihilate the entries in ${}^t\!\Gamma_{\mu^N}J_{n}$ of \thetag{\ref{eq:SchurQ_Pf_even}}
by elementary transformation of columns,
and annihilate the entries in $-J_{n}\Gamma_{\mu^N}$ of \thetag{\ref{eq:SchurQ_Pf_even}}
by elementary transformation of columns.
Then expand the Pfaffian $\Pf\left({\cal C}_{n,N}\right)$ along the first $N+1$ rows.
The case when $n$ is similar.
Perform the same operation on \thetag{\ref{eq:SchurQ_Pf_odd}}.
$\Box$
\end{demo}

%----------------------------------------------------------
% Proof of Corollary
%----------------------------------------------------------
\begin{demo}{Proof of Theorem~\ref{cor:sumQ_infinite}}
Perform the summations
\[
\sum_{0\leq k<l}\beta_{kl}
\det\begin{pmatrix}
x_{i}^{l}&x_{i}^{k}\\
x_{j}^{l}&x_{j}^{k}
\end{pmatrix}
\]
and
\[
\sum_{k=0}^{\infty}\beta_{-1,k}x_j^k,
\]
and apply Corollary~\ref{cor:sumQ_finite}.
The details are left to the reader
(cf. Proof of Theorem~2.1 in \cite{I}).
$\Box$
\end{demo}

To prove these theorems,
we need to cite a lemma from \cite{I}.
(See Corollary~3.3 of \cite{I} and Theorem~3.2 of \cite{IOTZ}.)
%----------------------------------------------------------
% Schur's Pfaffian
%----------------------------------------------------------
\begin{lemma}
\label{th:pfaff-det}
Let $n$ be a non-negative integer.
Let $X=(x_1,\dots,x_{2n})$, $A=(a_1,\dots,a_{2n})$, $B=(b_1,\dots,b_{2n})$, $C=(c_1,\dots,c_{2n})$ and $D=(d_1,\dots,d_{2n})$
be $2n$-tuples of variables.
Then
\begin{align}
&\Pf\left[\frac{(a_ib_j-a_jb_i)(c_id_j-c_jd_i)}
{(x_i-x_j)(1-tx_ix_j)}\right]_{1\leq i<j\leq 2n}
\nonumber\\
&\qquad\qquad=\frac{V^{n}(X,\pmb{1}+tX^2;A,B)V^{n}(X,\pmb{1}+tX^2;C,D)}
{\prod_{1\leq i<j\leq 2n}(x_i-x_j)(1-tx_ix_j)},
%\label{eq:second}
\end{align}
where $\pmb{1}+tX^2=(1+tx_1^2,\dots,1+tx_n^2)$.\par\noindent
In particular, we have
\begin{align}
\label{eq:cauchy}
&\Pf\left[\frac{a_ib_j-a_jb_i}
{1-tx_ix_j}\right]_{1\leq i<j\leq 2n}
=(-1)^{\binom{n}2}t^{\binom{n}2}\frac{V^{n}(X,\pmb{1}+tX^2;A,B)}{\prod_{1\leq i<j\leq 2n}(1-tx_ix_j)}.
\ \Box
\end{align}
\end{lemma}

%----------------------------------------------------------
% Proof of Theorem
%----------------------------------------------------------
\begin{demo}{Proof of Theorem~\ref{th:xi}}
First, assume $n$ is even.
Using the formula
\begin{equation}
\Pf(A+B)=\sum_{r=0}^{\lfloor n/2\rfloor}\sum_{I\in\binom{[n]}{2r}}(-1)^{|I|-r}
\Pf(A^{I}_{I})\Pf(B^{\overline I}_{\overline I}),
\label{eq:sum_Pf}
\end{equation}
where $\overline I$ denotes the complementary set of $I$,
we see that $\xi(a,b,c,d;X_{n})$ is equal to
\begin{equation*}
\sum_{r=0}^{\lfloor n/2\rfloor}
\sum_{I\in\binom{[n]}{2r}}(-1)^{|I|-r}
\prod_{{i,j\in I}\atop{i<j}}\frac{x_i+x_j}{x_i-x_j}
\Pf(\widetilde{v}_{ij})_{i,j\in I}.
\end{equation*}
Apply Lemma~\ref{th:pfaff-det} to obtain \thetag{\ref{eq:det_xi_even}}.
When $n$ is odd,
first expand the Pfaffian along the first row/column
and repeat the same argument.
$\Box$
\end{demo}

%----------------------------------------------------------
% Proof of Theorem
%----------------------------------------------------------
\begin{demo}{Proof of Theorem~\ref{th:zeta}}
%First, assume $n$ is even.
Note that the rank of the matrix $(u_{ij})_{1\leq i,j\leq n}$ is at most two.
Thus we have
\[
\Pf(u_{ij})_{1\leq i,j\leq n}
=\begin{cases}
\frac{a(x_1-x_2)\{1+b(x_1+x_2)+abx_1x_2\}}{(1-abx_1^2)(1-abx_2^2)}
&\text{ if $n=2$,}\\
0&\text{ otherwise.}
\end{cases}
\]
Using \thetag{\ref{eq:sum_Pf}},
we obtain
\begin{align*}
&\Pf\left(\gamma_{ij}\right)_{1\leq i,j\leq n}
=\Pf\left(\frac{x_{i}-x_{j}}{x_{i}+x_{j}}+v_{ij}z^2\right)_{1\leq i,j\leq n}\\
&+\sum_{1\leq k<l\leq n}(-1)^{k+l-1}
\frac{az(x_k-x_l)\{1+b(x_k+x_l)+abx_kx_l\}}{(1-abx_k^2)(1-abx_l^2)}
\Pf\left(\frac{x_{i}-x_{j}}{x_{i}+x_{j}}+v_{ij}z^2\right)_{{1\leq i,j\leq n}\atop{i,j\neq k,l}}.
\end{align*}
Use \thetag{\ref{eq:sum_Pf}} again to see that $\zeta(a,b,c,d;z;X_{n})$ is equal to
\begin{align*}
&\sum_{r=0}^{\lfloor n/2\rfloor}z^{2r}
\sum_{I\in\binom{[n]}{2r}}(-1)^{|I|-r}
\prod_{{i,j\in I}\atop{i<j}}\frac{x_i+x_j}{x_i-x_j}
\cdot\Pf(v_{ij})_{i,j\in I}\\
&+\sum_{1\leq k<l\leq n}(-1)^{k+l-1}
\frac{az(x_k-x_l)\{1+b(x_k+x_l)+abx_kx_l\}}{(1-abx_k^2)(1-abx_l^2)}\\
&\times\sum_{r=1}^{\lfloor n/2\rfloor}z^{2r-2}
\sum_{I'\in\binom{[n]-\{k,l\}}{2r-2}}
(-1)^{|I'|-r+1}
\prod_{{i,j\in I'}\atop{i<j}}\frac{x_i+x_j}{x_i-x_j}
\cdot\Pf(v_{ij})_{i,j\in I'}.
\end{align*}
Put $I=I'\cup\{k,l\}$ and apply Lemma~\ref{th:pfaff-det} to obtain \thetag{\ref{eq:det_zeta_even}}.
%The case when $n$ is odd can be proved similarly.
$\Box$
\end{demo}

\bigbreak
\noindent
{\bf Acknowledgment:}
The authors would like to express their gratitude to
Dr. Yasushi Kajihara
for his helpful comments and suggestions.

%%%%%%%%%%%%%%%%%%%%%%%%%%%%%%%%%%%%%%%%%%%%%%%%%%%%%%%%%%%%%%%%%%%%%%%%%%%%%
%
% References
%
%%%%%%%%%%%%%%%%%%%%%%%%%%%%%%%%%%%%%%%%%%%%%%%%%%%%%%%%%%%%%%%%%%%%%%%%%%%%%

%%%%%%%%%%%%%%%%%%%%%%%%%%%%%%%%%%%%%%%%%%%%%%%%%%%%%%%%%%%%%%%%%%%%%%%%%%%%%
%
% Reference
%
%%%%%%%%%%%%%%%%%%%%%%%%%%%%%%%%%%%%%%%%%%%%%%%%%%%%%%%%%%%%%%%%%%%%%%%%%%%%%

%%%%%%%%%%%%%%%%%%%%%%%%%%%%%%%%%%%%%%%%%%%%%%%%%%%%%%%%%%%%%%%%%%%%%%%%%%%%%
%
% Adrress
%
%%%%%%%%%%%%%%%%%%%%%%%%%%%%%%%%%%%%%%%%%%%%%%%%%%%%%%%%%%%%%%%%%%%%%%%%%%%%%

%\medskip
%\parindent=0mm
%
%Masao ISHIKAWA, Department of Mathematics, Faculty of Education, \\
%Tottori University, Tottori 680 8551, Japan
%
%E-mail: ishikawa@fed.tottori-u.ac.jp
%

%\medskip

%
%Jiang ZENG, Institut Camille Jordan,
%Universit\'e Claude Bernard Lyon I, \\
%43, boulevard du 11 novembre 1918
%69622 Villeurbanne Cedex, France
%
%E-mail: zeng@igd.univ-lyon1.fr

\end{document}